\documentclass[12pt,reqno]{amsart}

\usepackage{amsfonts}
\usepackage{amsmath}
\usepackage{amssymb}
\usepackage{amsthm}
\usepackage{bbm}
\usepackage{bm}
\usepackage{enumitem}
\usepackage{float}
\usepackage[T1]{fontenc} 
\usepackage[top=3cm, bottom=2.5cm, left=2.5cm, right=2.5cm]{geometry}
\usepackage{graphicx}
\usepackage[bookmarks=true,colorlinks=true,linkcolor=blue,citecolor=blue,urlcolor=blue,menucolor=blue,linktoc=all,pdfnewwindow=true]{hyperref}
\usepackage{mathrsfs}
\usepackage{mathtools}
\usepackage{physics}
\usepackage{relsize}
\usepackage[active]{srcltx}
\usepackage{thmtools}
\usepackage{tikz}
\usepackage{url}
\usepackage{xcolor}
\usepackage[all,cmtip]{xy}
\usepackage{xfrac}
\usepackage{hyperref}

\graphicspath{{Images/}}

\makeatletter
\def\th@plain{%
    \thm@notefont{}
    \itshape
}
\def\th@definition{%
    \thm@notefont{}
    \normalfont
}
\makeatother

\newtheorem{theorem}{Theorem}[section]
\newtheorem{lemma}[theorem]{Lemma}
\newtheorem{proposition}[theorem]{Proposition}

\newtheorem*{utheorem}{Theorem}

\theoremstyle{definition}

\theoremstyle{remark}
\newtheorem{remark}[theorem]{Remark}

\newcommand{\M}{\mathbf{M}}
\newcommand{\GL}{\mathbf{GL}}
\newcommand{\SL}{\mathbf{SL}}

\renewcommand{\d}{\, \mathrm{d}}

\def\veca{{\text{\boldmath$a$}}}
\def\vecb{{\text{\boldmath$b$}}}
\def\vece{{\text{\boldmath$e$}}}
\def\vecf{{\text{\boldmath$f$}}}

\def\vecv{{\text{\boldmath$v$}}}
\def\vecw{{\text{\boldmath$w$}}}
\def\vecx{{\text{\boldmath$x$}}}
\def\vecy{{\text{\boldmath$y$}}}

\newcommand{\bn}{\mathbf{0}}
\newcommand{\ve}{\varepsilon}

\newcommand{\frakI}{\mathfrak{I}}
\newcommand{\frakp}{\mathfrak{p}}
\newcommand{\bG}{\mathbf{G}}
\newcommand{\bH}{\mathbf{H}}
\newcommand{\bM}{\mathbf{M}}
\newcommand{\bU}{\mathbf{U}}

\newcommand{\td}{\widetilde{d}}
\newcommand{\tD}{\widetilde{D}}
\newcommand{\tx}{\widetilde{x}}
\newcommand{\ty}{\widetilde{y}}
\newcommand{\A}{\mathbb{A}}
\newcommand{\CC}{\mathbb{C}}
\newcommand{\R}{\mathbb{R}}
\newcommand{\Z}{\mathbb{Z}}
\newcommand{\Q}{\mathbb{Q}}
\newcommand{\matr}[4]{\left( \begin{matrix} #1 & #2 \\ #3 & #4 \end{matrix} \right) }
\newcommand{\cmatr}[2]{\left( \begin{matrix} #1 \\ #2 \end{matrix} \right) }
\newcommand{\col}{\: : \:}
\newcommand{\hx}{\widehat{x}}

\def\scrF{{\mathcal F}}
\def\scrO{{\mathcal O}}
\def\scrP{{\mathcal P}}
\def\scrX{{\mathcal X}}
\def\scrY{{\mathcal Y}}
\def\new{\operatorname{new}}
\DeclareMathOperator{\C}{C}

\DeclareMathOperator{\vol}{vol}
\DeclareMathOperator*{\supp}{supp}

\DeclareMathOperator{\ord}{ord}

\numberwithin{equation}{section}

\begin{document}


\title[]{Adelic Rogers' formula and an application}

\author{Mahbub Alam}

\address{\parbox{\linewidth}{\textbf{Mahbub Alam}
    Department of Mathematics,
    Uppsala University, Sweden\\[-4pt]
\url{https://sites.google.com/view/mahbubweb}\\[-10pt]}}

\email{mahbub.dta@gmail.com, mahbub.alam@math.uu.se}

\author{Andreas Str\"ombergsson}

\address{\textbf{Andreas Str\"ombergsson} \\
    Department of Mathematics,
    Uppsala University, Sweden}
\email{astrombe@math.uu.se}

\date{}

\begin{abstract}
    Recently, Seungki Kim \cite{kim24} proved an extension of 
Rogers' mean value formula \cite{rogers55a} 
to the adeles of an arbitrary number field.
In this paper we give a new proof Kim's formula,
and give a criterion ensuring convergence in this formula.
We also discuss one application, namely diophantine approximation 
over imaginary quadratic number fields
with congruence conditions, where we prove an analogue of a famous counting result of \mbox{W.\ M.\ Schimdt}.
\end{abstract}

\thanks{Both Alam and Str\"ombergsson were supported by the Knut and Alice Wallenberg Foundation}

\maketitle

\setcounter{tocdepth}{1}
\tableofcontents

\section{Introduction}\label{sec:intro}

A fundamental integration formula due to Siegel \cite{siegel45}
states that for any measurable function $f:\R^d\to\R_{\geq0}$ ($d\geq2$),
\begin{align}\label{Siegelformula}
\int_{\scrY_d}\sum_{\vecv\in L\smallsetminus\{\bn\}} f(\vecv)\d\mu(L)  
=\int_{\R^d}f(\vecx)\d\vecx.
\end{align}
Here $\scrY_d$ is the space of lattices in $\R^d$ of covolume one,
equipped with its invariant probability measure $\mu$,
and $\d\vecx$ denotes Lebesgue measure on $\R^d$.
The integrand in the left hand side of \eqref{Siegelformula},
i.e., the function on $\scrY_d$ mapping $L$ 
to $\sum_{\vecv\in L\smallsetminus\{\bn\}} f(\vecv)$,
is often called the Siegel transform of $f$.

Rogers \cite{rogers55a}\footnote{See also the papers
\cite{wS57,wS58a}
and
\cite{macbeathrogers58}
for alternative and corrected proofs.}
proved a generalization of \eqref{Siegelformula}
which for any $1\leq k<d$ 
gives an explicit formula for 
the average of a $k$-fold 
sum over $L$, that is, a formula for
\begin{align}\label{Rogersformula}
\int_{\scrY_d}\sum_{\vecv_1,\ldots,\vecv_k\in L\smallsetminus \{\bn\}}f(\vecv_1,\ldots,\vecv_k)\,d\mu(L),
\end{align}
for any $f:(\R^d)^k\to\R_{\geq0}$.
(The formula in \cite{rogers55a} has 
$\sum_{\vecv_1,\ldots,\vecv_k\in L}$ instead of 
$\sum_{\vecv_1,\ldots,\vecv_k\in L\smallsetminus \{\bn\}}$ inside the integral; 
however it is easy to translate between the two.)
For comparison with formulas appearing below, 
recall that $\scrY_d$ can be
identified with the homogeneous space $\SL_d(\R)/\SL_d(\Z)$
through the map $g\,\SL_d(\Z)\mapsto g\Z^d$.
In this way, identifying also $(\R^d)^k$ with the space of real $d\times k$ matrices,
$\bM_{d,k}(\R)$, the integral in \eqref{Rogersformula} can be written:
\begin{align}
\int_{\scrY_d}\sum_{\substack{X \in \bM_{d,k}(\Z) \\ \text{columns nonzero}}} f(gX)\,d\mu(g).
\end{align}

The integration formulas by Siegel and Rogers
have found a number of applications over the years,
to lattice problems and the geometry of numbers, 
and to problems in Diophantine approximation;
see, e.g.,
\cite{rogers56},
\cite{schmidt58b},
\cite{pSaS2006},
\cite{aS2011o},
\cite{jathreyagmargulis09},
\cite{bjorklundgorodnik19},
\cite{strombergssonsodergren19}.
In recent years, many variants of Rogers' formula have been developed,
with several new applications;
see, e.g.,
\cite{dKsY2019},
\cite{ghoshkelmeryu2020},
\cite{han22},
\cite{alamghoshhan21},
\cite{alamghoshyu21},
\cite{aSaS2022},
\cite{kim24},
\cite{nH2023},
\cite{gargavasv23},
\cite{nG2023c}.

In particular, in \cite{kim24},
Seungki Kim proved an extension of Rogers' integration formula to the adeles of an arbitrary number field.
One of our main purposes 
is to present an alternative proof of Kim's formula.
We start by recalling its statement.

\subsection{Adelic Rogers' formula}\label{ssec:arf}

Let $F$ be a number field,
let $\A_F$ be the adele ring of $F$,
and set
\begin{align*}
G_d = \{g \in \GL_d(\mathbb{A}_{F}) : \norm{\det g}_{\mathbb{A}_F} = 1\}.
\end{align*}
Then $\GL_d(F)$ is a lattice in $G_d$.   
We let $\mathcal{X}_d$ be the homogeneous space $G_d/\GL_d(F)$,
and let $\mu_d^1$ be the (left and right) Haar measure on $G_d$
normalized so that $\mu_d^1(\scrX_d)=1$.
Denote by $\alpha_F$ the Tamagawa measure on $\mathbb{A}_F$,
i.e., the Haar measure on $\langle\A_F,+\rangle$
normalized so that $\alpha_F(\A_F/F)=1$.
For any positive integers $d,k$, let $\bM_{d,k}(\A_F)$ be the group of 
$d\times k$ matrices with entries in $\A_F$,
and let $\alpha_F^{dk}$ be the Haar measure on $\bM_{d,k}(\A_F)$
which is the $dk$-fold product measure of $\alpha_F$.

The following is a restatement of Kim's \cite[Theorem 1.2]{kim24},
with the addition of a statement about the convergence of the integrals involved. 
\begin{theorem}\label{thm:arf}
Let $1\leq k<d$, and let $f : \bM_{d,k}(\A_F)\to\R_{\geq0}$ be a Borel measurable function.
Then
\begin{align}\label{eq:arf}
            \int_{\mathcal{X}_d} \sum_{\substack{X \in \bM_{d,k}(F) \\ \text{columns nonzero}}} f(gX) \d{\mu_d^1(g)} 
= \int_{\bM_{d,k}(\A_F)} f(X) \d{\alpha^{dk}_{F}(X)}
+\sum_{m=1}^{k-1}\sum_D\int_{\bM_{d,m}(\A_F)} f(XD) \d{\alpha^{dm}_{F}(X)},
\end{align}
    where the sum over $D$ is over all $m \times k$ row-reduced echelon forms over $F$ of rank $m$ with no zero columns.
The relation in \eqref{eq:arf} is an equality of numbers in $\R_{\geq0}\cup\{+\infty\}$,
i.e., either both sides are finite and equal or else both are $+\infty$.
Furthermore, if $f$ is bounded and has compact support, 
then both sides 
in \eqref{eq:arf} are finite.
\end{theorem}
\begin{remark}\label{thm:arfcorrREM}
We have corrected the statement of \cite[Theorem 1.2]{kim24}
by requiring $D$ in \eqref{eq:arf} to have \textit{no zero columns.}
(See \autoref{KIMcomparisonREM1} below.)
\end{remark}
\begin{remark}
As an immediate consequence of \autoref{thm:arf},
for any $k,d,f$ as in that theorem we also have
\begin{equation}\label{eq:arf0}
         \int_{\mathcal{X}_d} \sum_{\substack{X \in \bM_{d,k}(F) \\ \text{columns indep}}} f(gX) \d{\mu_d^1(g)} 
= \int_{\bM_{d,k}(\A_F)} f(X) \d{\alpha^{dk}_{F}(X)},
\end{equation}
where the sum over $X$ is over all $X\in \bM_{d,k}(F)$
whose columns are linearly independent over $F$.
In Kim's paper, \eqref{eq:arf0} was stated in the first part of
\cite[Theorem 1.2]{kim24}.
To prove \eqref{eq:arf0},
for given $f$ we 
define $f_1 : \bM_{d,k}(\A_F)\to\R_{\geq0}$
through $f_1(X)=f(X)$ if the columns of $X$ are linearly independent over $F$,
otherwise $f_1(X)=0$.
Then \eqref{eq:arf} applied to $f_1$ gives the relation in
\eqref{eq:arf0}.
\end{remark}
\begin{remark}\label{thm:arfposnegREM}
From the convergence statement in the last part of
\autoref{thm:arf} it follows that
for bounded functions $f$ of compact support,
we may remove the assumption that $f$ is non-negative,
i.e., the formula \eqref{eq:arf}
holds for any complex-valued Borel measurable function
$f$ on $\bM_{d,k}(\A_F)$ which is bounded and of compact support,
with both sides of \eqref{eq:arf} being nicely absolutely convergent.
\end{remark}

\begin{remark}
For the original formula by Rogers, \cite[Theorem 4]{rogers55a},
the analogue of the convergence statement in the last part of \autoref{thm:arf}
was proved by W.\ M.\ Schmidt \cite{schmidt58a}.
Also for the (non-adelic) analogue of Rogers' formula over an arbitrary number field,
which was obtained recently by Nathan Hughes \cite{nH2023},
the analogous convergence statement has been proved, in
Gargava-Serban-Viazovska
\cite[Cor.\ 11]{gargavasv23}.
A crucial ingredient,
both in our proof of convergence and in that of 
\cite[Cor.\ 11]{gargavasv23},
is an estimate due to W.\ M.\ Schmidt, \cite[Theorem 3]{schmidt67},
on the number of linear subspaces of $F^k$ of bounded height.
\end{remark}

\subsection{Diophantine approximation on $\mathbb{C}$ with congruence conditions}\label{ssec:cdioc}
Diophantine approximation on $\mathbb{C}$ has been of interest at least since the late 1800's in the works of Hurwitz \cite{hurwitz1887} (who discusses continued fractions of complex numbers with Gaussian and Eisenstein integers) and in the works of Ford \cite{ford18, ford25}.
 Over the years, diophantine approximation on $\mathbb{C}$ via continued fractions of Gaussian integers was discussed by, among others, W.\ J.\ LeVeque \cite{leveque52} and A.\ L.\ Schmidt \cite{aschmidt75}.
In recent years, diophantine approximation over number fields (which includes the above) and related topics were studied by Qu\^{e}me \cite{rqueme89}, Ly \cite{ly16} and Kleinbock-Ly \cite{kleinbockly16}, Kleinbock-Shi-Tomanov \cite{kleinbockshitomanov17} and Alam-Ghosh \cite{alamghosh20} among others.
Before discussing known results for motivation we introduce the problem that we consider in the present paper.

We treat the problem of diophantine approximation on $\mathbb{C}$ with congruence conditions.
More precisely, we consider imaginary quadratic number fields $F := \mathbb{Q}(\sqrt{-D})$, for a square-free positive integer $D$, and use its ring of integers, denoted by $\mathcal{O}_F = \mathcal{O}$, for diophantine approximation on $\mathbb{C}$.
Note that $\mathbb{C}$ is the completion of $F$ under the (archimedean) valuation defined on $F$.
For $\vartheta \in \M_{m \times n}(\mathbb{C})$, a decreasing function
$\psi:[1,\infty)\to(0,\infty)$
and with $(\bm{p}, \bm{q}) \in \mathcal{O}^{m} \times \mathcal{O}^{n}$, consider
\begin{equation}\label{eq:cdioc0}
    \begin{gathered}
        \norm{\vartheta\bm{q} + \bm{p}}^{m}_\infty \leq \psi(\norm{\bm{q}}_{\infty}^{n}).
    \end{gathered}
\end{equation}
Recall that $\norm{ \cdot }_\infty$ denotes the absolute value in $\mathbb{C}$, i.e., the square of the usual distance function on $\mathbb{C}$, and by abuse of notation, let it also denote the `sup norm' on $\mathbb{C}^{k}$, for any $k > 0$, i.e., $\norm{\bm{x}}_\infty := \sup_{i} \norm{x_i}_\infty$ for $\bm{x} \in \mathbb{C}^{k}$.
Finding a criterion in terms of $\psi$ that ensures infinitely many solutions $(\bm{p}, \bm{q})$ to \eqref{eq:cdioc0} for almost every $\vartheta \in \M_{m \times n}(\mathbb{C})$ is one of the fundamental problems in diophantine approximation.
The following result is known in this direction:
\begin{utheorem}[\cite{leveque52, sullivan82, ly16}]\label{thm:cdioc0}
    \begin{enumerate}[label=(\alph*)]
        \item If $\int_{1}^{\infty} \psi(t) \d{t} < \infty$, then for (Lebesgue) almost every $\vartheta \in \M_{m \times n}(\mathbb{C})$, 
\eqref{eq:cdioc0} has at most finitely many solutions.
        \item If $\int_{1}^{\infty} \psi(t) \d{t} = \infty$, then for almost every $\vartheta \in \M_{m \times n}(\mathbb{C})$,
\eqref{eq:cdioc0} has infinitely many solutions.
    \end{enumerate}
\end{utheorem}
Our main result is towards a quantitative generalization of the above theorem with additional congruence conditions.
For $T > 1$, a fixed $\bm{v} \in \mathcal{O}^{m} \times \mathcal{O}^{n} = \mathcal{O}^{d}$ and an ideal $\mathfrak{I}$ of $\mathcal{O}$, let us denote by $\mathcal{N}(\vartheta, T, \bm{v}, \mathfrak{I})$ the number of solutions 
$(\bm{p}, \bm{q})\in \mathcal{O}^{m} \times \mathcal{O}^{n}$ to
\begin{equation}\label{eq:cdioc}
    \begin{gathered}
        \norm{\vartheta\bm{q} + \bm{p}}^{m}_\infty \leq \psi(\norm{\bm{q}}_{\infty}^{n}), \\
        1\leq \norm{\bm{q}}_{\infty}^{n} < T, \\
        (\bm{p}, \bm{q}) = \bm{v} \pmod {\mathfrak{I}}.
    \end{gathered}
\end{equation}

We let $d = m + n$ and define
\[
    E_{T, \infty} := \{(\bm{x}, \bm{y}) \in \mathbb{C}^{m} \times \mathbb{C}^{n} = \mathbb{C}^{d}: \norm{\bm{x}}_\infty^m \leq \psi(\norm{\bm{y}}_\infty^n) \ \text{and} \ 1 \leq \norm{\bm{y}}_\infty^n < T\}.
\]
We will show in \S\ref{ssec:dcx} that
\[
    \vol(E_{T, \infty}) = \pi^d \int_{1}^{T} \psi(t) \d{t},
\]
where ``$\vol$'' is the standard $2d$-dimensional volume measure on $\mathbb{C}^d=(\R\oplus i\R)^d$.
Our main result relates the set $E_{T, \infty}$ and its volume to the solutions of \eqref{eq:cdioc}:
\begin{theorem}\label{thm:cdioc}
Let $m,n\geq1$ and assume that $d:=m+n\geq3$.
Let $F$ be an imaginary quadratic number field and $\scrO$ its ring of integers;
let $\frakI$ be a non-zero ideal of $\scrO$, and let $\vecv\in\scrO^m\times\scrO^n$.
Let $\psi:[1,\infty)\to(0,\infty)$ be a decreasing function
satisfying $\int_{1}^{\infty} \psi(t) \d{t} = \infty$. 
Then for Lebesgue almost every $\vartheta \in \M_{m \times n}(\mathbb{C})$,
\begin{align}\label{thm:cdiocRES}
\mathcal{N}(\vartheta, T, \bm{v}, \mathfrak{I}) \thicksim 2^d |\Delta_F|^{-d/2} N(\mathfrak{I})^{-d}\vol(E_{T, \infty}), \qquad \text{as} \ T \to \infty.
\end{align}
Here $\Delta_F$ is the discriminant of $F$,
and $N(\mathfrak{I})=\#(\scrO/\frakI)$ is the norm of the ideal $\mathfrak{I}$.
\end{theorem}
(We adopt the standard convention that $f(T) \thicksim g(T)$ as $T \to \infty$ means that $\frac{f(T)}{g(T)} \to 1$ as $T \to \infty$.)

\vspace{5pt}

We conclude this introductory section by discussing previously known results related to our \autoref{thm:cdioc}.
We will here consider the cases of $\R$ and $\CC$ in parallel.
In diophantine approximation on $\M_{m \times n}(\mathbb{R})$,
one is interested in the 
solutions $(\bm{p}, \bm{q}) \in \mathbb{Z}^{m} \times \mathbb{Z}^{n}$ to the inequality
\begin{equation}\label{eq:drmn}
    \norm{\vartheta\bm{q} + \bm{p}}^m \leq \psi(\norm{\bm{q}}^n), \\
\end{equation}
where $\vartheta \in \M_{m \times n}(\mathbb{R})$
and $\psi:[1,\infty)\to(0,\infty)$ is a decreasing function; 
and in quantitative diophantine approximation one is interested in the number of integral solutions satisfying
\begin{equation}\label{eq:drmnq}
    \begin{gathered}
        \norm{\vartheta\bm{q} + \bm{p}}^m \leq \psi(\norm{\bm{q}}^n), \\
        1 \leq \norm{\bm{q}}^n < T,
    \end{gathered}
\end{equation}
for $T > 1$.
Here $\norm{ \cdot }$ may denote any fixed chosen norm (on $\mathbb{R}^{n}$ and $\mathbb{R}^{m}$), but for simplicity, we take it to be the sup norm on both spaces.
Analogously, for the complex case, with $\vartheta \in \M_{m \times n}(\mathbb{C})$ and $(\bm{p}, \bm{q}) \in \mathcal{O}^{m} \times \mathcal{O}^{n}$, one considers
\begin{equation}\label{eq:dioc}
    \begin{gathered}
        \norm{\vartheta\bm{q} + \bm{p}}^{m}_\infty \leq \psi(\norm{\bm{q}}_{\infty}^{n}),
    \end{gathered}
\end{equation}
and the quantitative case
\begin{equation}\label{eq:diocq}
    \begin{gathered}
        \norm{\vartheta\bm{q} + \bm{p}}^{m}_\infty \leq \psi(\norm{\bm{q}}_{\infty}^{n}), \\
        1\leq \norm{\bm{q}}_{\infty}^{n} < T.
    \end{gathered}
\end{equation}

Khinchin-Groshev \cite{khinchin26, groshev38} showed that \eqref{eq:drmn} has infinitely many (resp.\ finitely many) solutions for almost every $\vartheta \in \M_{m \times n}(\mathbb{R})$ if and only if
\begin{equation}\label{eq:psi}
    \int_{1}^{\infty} \psi(t) \d{t}
\end{equation}
diverges (resp.\ converges).
For diophantine approximation on $\mathbb{C}$ with the ring of integers of $\mathbb{Q}(\sqrt{-D})$, Sullivan \cite[Theorem 1]{sullivan82} used the notion of disjoint spheres and geodesic excursions to cuspidal neighborhoods of a finite volume quotient of hyperbolic space to derive a Khinchin-type theorem relating the divergence (resp.\ convergence) of \eqref{eq:psi} to the existence of infinitely many (resp.\ finitely many) solutions of \eqref{eq:dioc} for almost every $\vartheta \in \mathbb{C}$.
For the specific case of $D = 1$, a Khinchin-type theorem was derived earlier by LeVeque \cite{leveque52}.
(To be precise, both of them considered the case $m = n = 1$.)
For general $m$ and $n$, a Khinchin-type theorem has been proved in \cite{ly16}, in fact for general number fields.
In the light of these Khinchin-type theorems,
the assumption in \autoref{thm:cdioc} is necessary if we want $\mathcal{N}(\vartheta, T, \bm{v}, \mathfrak{I})$ to be unbounded as function of $T$.

W.\ M.\ Schmidt \cite[Theorem 1]{schmidt60b} strengthened 
the Khinchin-Groshev theorem into the following quantitative statement
(here stated with a slightly weaker error term than in \cite{schmidt60b}):  
The number of solutions to \eqref{eq:drmnq} is 
\[
    2^{d} \Psi(T) + O(\Psi(T)^{1/2 + \varepsilon}) \: \text{ for every} \ \varepsilon > 0,
\]
where
\[
    \Psi(T) = \int_{1}^{T} \psi(t) \d{t}.
\]

By applying LeVeque and Sullivan's theorem for $\psi(t) = \frac{c}{t^s}$, 
we see that $s = 1$ is the critical exponent for which \eqref{eq:dioc} can have infinitely many solutions for almost every $\vartheta$.
For $m = n = 1$ and this critical exponent, i.e., for $\psi(t) = \frac{c}{t}$, Nakada \cite{nakada88} proved that the number of solutions to \eqref{eq:diocq} is asymptotic to $2^dcT$.
In other words, he proved an analogue of W.\ M.\ Schmidt's result, without the error term, for the specific case when $\psi(t) = \frac{c}{t}$.
\cite{alamghosh20} extended Nakada's theorem to any number field (still with $\psi(t) = \frac{c}{t}$) by proving an equidistribution result on a space of unimodular lattices associated to the number field.
%

To the best of our knowledge, an analogue of W.\ M.\ Schmidt's result for 
number fields other than $\Q$ 
with general $\psi$ is unknown.
Our result \autoref{thm:cdioc} gives an analogue of \mbox{W.\ M.\ Schmidt's} result for the case of an imaginary quadratic number field
with general $\psi$ but without an error term (with $\mathfrak{I} = \mathcal{O}$).


\textbf{Congruence conditions:}
A refinement of the above problem \eqref{eq:drmnq} is the inquiry into the various congruence relations that the approximates $(\bm{p}, \bm{q})$ might satisfy.
For the real case, one might guess that for any fixed positive integer $N$, the solutions $(\bm{p}, \bm{q})$ to \eqref{eq:drmnq} would be equally distributed in all congruence classes modulo $N$, as $T \to \infty$.
For fixed $\bm{v} \in \mathbb{Z}^{m} \times \mathbb{Z}^{n} = \mathbb{Z}^{d}$ and $N \in \mathbb{Z}_{>0}$ we consider:
\begin{equation}\label{eq:dc}
    \begin{gathered}
        \norm{\vartheta\bm{q} + \bm{p}}^m \leq \psi(\norm{\bm{q}}^n), \\
        1 \leq \norm{\bm{q}}^n < T, \\
        (\bm{p}, \bm{q}) = \bm{v} \pmod N.
    \end{gathered}
\end{equation}
Analogously, for the complex case, we consider \eqref{eq:cdioc}.

The real case with $m = n = 1$ was studied by Sz\"usz \cite{szusz62} where he counts the number of solutions to \eqref{eq:dc}.
Recently a Khinchin-Groshev-type theorem (with general $m, n$) in this set-up was proved in \cite{nesharimruhrshi2020} using dynamics on a space of affine lattices in $\mathbb{R}^{d}$.
Moreover, in \cite{alamghoshyu21} the problem of counting the number of solutions to \eqref{eq:dc} was translated into an affine lattice point counting problem and a Schmidt-type theorem without the error term was proved.
More precisely, they consider the affine lattices of the form ${g(\mathbb{Z}^{d} + \frac{\bm{v}}{N})}$ where $g \in \SL_d( \mathbb{R})$.
The collection $\mathcal{Y}_{N, \bm{v}}$ of such affine lattices can be identified with $\SL_d( \mathbb{R})/\Gamma_{N, \bm{v}}$, where $\Gamma_{N, \bm{v}} := \{\gamma \in \SL_d( \mathbb{Z}) : \gamma\bm{v} = \bm{v} \pmod N\}$ is the stabilizer of the affine lattice $\mathbb{Z}^{d} + \frac{\bm{v}}{N}$ in $\SL_{d}(\mathbb{R})$.
It can then be easily shown (\cite[\S2.3]{alamghoshyu21}) that the number, $\mathcal{N}_\mathbb{R}(\vartheta, T, \bm{v}, N)$, of solutions to \eqref{eq:dc} is
\[
    \#{\left(u(\vartheta){\left(\mathbb{Z}^{d} + \frac{\bm{v}}{N}\right)} \cap N^{-1}E_T\right)},
\]
and they proved that this number is related to $\vol(E_T) = 2^d \Psi(T)$.
\begin{theorem}[\cite{alamghoshyu21}, Theorem 1.1]\label{thm:dioccagy}
    For almost every $\vartheta \in \M_{m \times n}(\mathbb{R})$
    \[
        \mathcal{N}_\mathbb{R}(\vartheta, T, \bm{v}, N) = \#{\left(u(\vartheta){\left(\mathbb{Z}^{d} + \frac{\bm{v}}{N}\right)} \cap N^{-1}E_T\right)} \thicksim N^{-d} \vol(E_T).
    \]
\end{theorem}
\vspace{2em}



\subsection{Notation}\label{notationSEC}
We will mainly use the same notation as in \cite{kim24}.
In particular, we write $\scrP_F$ for the set of all places of $F$,
so that $\A_F$ is the restricted product $\prod_{\nu\in\scrP_F}' F_\nu$.
We also write $\mathbf{f}\subseteq\scrP_F$ for the subset of finite places,
and set $\scrP_{F,\infty}=\scrP_F\smallsetminus\mathbf{f}$.
For each $\nu\in\mathbf{f}$ we let $\scrO_\nu$ be the maximal compact subring of $F_\nu$,
and choose a uniformizer $\pi_\nu\in\scrO_\nu$, so that
$(\pi_\nu)$ is the unique maximal ideal of $\scrO_\nu$.


For each place $\nu\in\scrP_F$ we define $\alpha_\nu$
to be the Haar measure on $\langle F_\nu,+\rangle$,
normalized so that $\alpha_\nu(\scrO_\nu)=1$
if $\nu$ is finite;
to be the standard Lebesgue measure on $\R$ if $\nu$ is real;
and to be twice the standard Lebesgue measure on $\mathbb{C}=\R\oplus i\R$
if $\nu$ is complex.
The Tamagawa measure $\alpha_F$ of $\A_F$ is then given by
\begin{align}\label{alphaFDEF}
\alpha_F=|\Delta_F|^{-\frac12}\prod_{\nu\in\scrP_F}\alpha_\nu,
\end{align}
where $\Delta_F$ is the discriminant of $F$.
We recall that $\alpha_F(\A_F/F)=1$.

Also for each place $\nu\in\scrP_F$ we define the absolute value
$\|\cdot\|_\nu:F_\nu\to\R_{\geq0}$ so that 
$\alpha_\nu(xA)=\|x\|_\nu\, \alpha_\nu(A)$ for every $x\in F_\nu$
and every Borel set $A\subseteq F_\nu$.
Note that we then have
\begin{align*}
\|x\|_{\A_F}=\prod_{\nu\in\scrP_F}\|x_\nu\|_\nu
\quad\text{for all }\: x=(x_\nu)\in\A_F^\times.   
\end{align*}

We denote by $\scrO_F$ the ring of integers of $F$.


\section{New proof of Adelic Rogers' formula}
\label{newproofSEC}

We will prove \autoref{thm:arf} using a technique that is due to Weil;
see \cite{weil82}.
The same technique can also be used to prove the original formula by Rogers,
\cite[Theorem 4]{rogers55a}.

Let $F$ be a number field.
Let $1\leq k<d$.

Let $\bG=\GL_d$ as an algebraic group defined over $F$,
and let $\bM_{d,k}$ be affine $dk$-space over $F$,
which we view as the space of $d\times k$ matrices
and on which $\bG$ acts by 
left multiplication.
Let $p_0$ be the point $\cmatr{I_k}0$ in $\bM_{d,k}$
(block matrix notation),
and let $\bH$ be the stabilizer of $p_0$ in $\bG$;
this is an algebraic subgroup of $\bG$ defined over $F$.
The $\bG$-orbit of $p_0$ equals 
the open subvariety $\bU$ of $\bM_{d,k}$
consisting of all matrices in $\bM_{d,k}$ of full rank,
and by \cite[Prop.\ 6.7]{aB91},
the morphism $\pi:\bG\to\bU$, $\pi(g)=gp_0$, 
realizes the quotient $\bG/\bH$;
we will from now on identify $\bG/\bH$ with $\bU$.

Let $\bG(\A_F)$, $\bH(\A_F)$ and $\bU(\A_F)$ be the adele-spaces attached to
$\bG$, $\bH$ and $\bU$ 
over $F$ \cite[Ch.\ 1.2]{weil82}.
Thus $\bG(\A_F)=\GL_d(\A_F)$.
In order to describe $\bU(\A_F)$ explicitly,
let $B$ be the family of subsets $\beta\subseteq\{1,\ldots,d\}$
of cardinality $k$, and for each $\beta\in B$ and $x\in\bM_{d,k}$ 
(or $x\in\bM_{d,k}(R)$ for any commutative ring $R$),
let us write $x_\beta$ for the determinant of the $k\times k$ submatrix of $x$
consisting of the rows with indices $i\in\beta$;
then $x\in\bM_{d,k}$ lies in $\bU$ if and only if $x_\beta\neq0$ 
for at least one $\beta\in B$.
Next, for each $\nu\in\mathbf{f}$ we set 
\begin{align*}
C_{\bU,\nu}:=\{x\in\bM_{d,k}(\scrO_\nu)\col[\exists\beta\in B\text{ s.t. } x_\beta\in\scrO_\nu^\times]\};
\end{align*}
this is a compact and open subset of $\bU(F_\nu)$.
For each finite subset $S\subseteq\scrP_F$ containing all infinite places,
we set
\begin{align*}
\bU_S:=\prod_{\nu\in S}\bU(F_\nu)\times \prod_{\nu\in\scrP_F\smallsetminus S} C_{\bU,\nu},
\end{align*}
and equip $\bU_S$ with the product topology.
Note that each set $\bU_S$ is a subset of $\bM_{d,k}(\A_F)$.
Now the construction of adele-space in 
\cite[Ch.\ 1.2]{weil82} gives that
$\bU(\A_F)$ equals the union of all these set $\bU_S$,
and is equipped with the inductive limit topology.
Thus a set $V\subseteq\bU(\A_F)$ is open if and only if $V\cap\bU_S$ is open in $\bU_S$
for every $S$.

One verifies that the map $g\mapsto gp_0$ from $\bG(\A_F)$ to $\bU(\A_F)$ induces 
a homeomorphism from $\bG(\A_F)/\bH(\A_F)$ onto $\bU(\A_F)$
(this also follows from \cite[Thm.\ 2.4.2]{weil82}).

Next we introduce a suitably normalized Haar measure $\mu_d$ on $\GL_d(\A_F)$,
following \cite[Ch.\ I.4]{jO84}.
For each $\nu\in\scrP_F$, let $\mu_{d,\nu}$ be the 
Haar measure on
$\GL_d(F_\nu)$ given by
\begin{align*}
d\mu_{d,\nu}(g)=\|\det(g)\|_{\nu}^{-d}\, \d \alpha_\nu(g_{1,1})\d \alpha_\nu(g_{1,2})\,\cdots\d \alpha_\nu(g_{n,n}),
\end{align*}
where $g_{i,j}$ are the entries of the matrix $g$.
Then define the (left and right) Haar measure $\mu_d$ on
$\GL_d(\A_F)$ by
\begin{align}\label{munDEF}
\mu_d=\bigl|\Delta_F\bigr|^{-d^2/2}\rho_F^{-1}
\prod_{\nu\mid\infty}\mu_{d,\nu}\times\prod_{\nu\in\mathbf{f}}(1-\|\pi_\nu\|_\nu)^{-1}\mu_{d,\nu},
\end{align}
where $\rho_F$ is the residue of the Dedekind zeta function $\zeta_F(s)$ at $s=1$.
To be precise, for each finite subset $S\subseteq\scrP_F$ containing all infinite places, 
the restriction of $\mu_d$ to
$\prod_{\nu\in S}\bG(F_\nu)\times\prod_{\nu\in\scrP_F\smallsetminus S}\bG(\scrO_\nu)$
equals the product measure in the right hand side of \eqref{munDEF}.
Note that $\mu_d$ equals the measure $\mu_{\GL_d}$
in \cite[Def.\ 4.7]{jO84}
(indeed, apply \cite[4.6--7]{jO84}
with $S=\scrP_F\smallsetminus\mathbf{f}$
and $\omega$ equal to the differential form $\omega_d$ in 
the proof of \autoref{GAFquotLEM2} below).

Next, as in 
\cite[Ch.\ I.5]{jO84},
we set
\begin{align*}
\bG^1(\A_F)=\{g\in\bG(\A_F)\col \|\det(g)\|_{\A_F}=1\}=G_d,
\end{align*}
and let $\mu_d^1$ be the Haar measure on $G_d$
normalized so that for all $f\in\C_c(\bG(\A_F))$,
\begin{align}\label{omegan1DEF}
\int_{\bG(\A_F)}f(g)\d \mu_d(g)
=\int_{\R^+}\int_{G_d} f(g\delta_t)\d \mu_d^1(g)\,\frac{\d t}t,
\end{align}
where for each $t\in\R^+$, $\delta_t$ is an element of
$\bG(\A_F)$ satisfying $\|\det(\delta_t)\|_{\A_F}=t$.
Note that the inner integral in \eqref{omegan1DEF},
$\int_{G_d} f(g\delta_t)\d \mu_d^1(g)$, only depends on 
$t$ and $\mu_d^1$, and not on the (continuous) choice of $t \mapsto \delta_t$.

It is known that the Tamagawa number of $\SL_d$ equals $1$
\cite[Theorem 3.3.1]{weil82},
and also the Tamagawa number of the multiplicative group
$\mathbb{G}_m=\GL_1$    
equals 1 \cite[Ex.\ 5.14(2)]{jO84}, \cite[Ch.\ VII.5--6]{weil95}.
Hence by \cite[Th\'eor\`eme 5.3]{jO84}
applied to the exact sequence
\begin{align*}
1\to\SL_d\to\GL_d\stackrel{\det}{\longrightarrow}\mathbb{G}_m\to1,
\end{align*}
we have $\tau_{\GL_d}=\tau_{\SL_d}\cdot\tau_{\mathbb{G}_m}=1$,
i.e., also the Tamagawa number of $\GL_d$ equals $1$.
By \cite[Def.\ 5.12]{jO84} and since $\bG^1(\A_F)=G_d$, this means that
\begin{align}\label{TamagawaGLneq1}
\mu_d^1(G_d/\GL_d(F))=1.
\end{align}

Next we set
\begin{align*}
\bH^1(\A_F)=\bH(\A_F)\cap G_d   
=\left\{\matr{I_k}Z0h\col Z\in\bM_{k,d-k}(\A_F),\: h\in G_{d-k}\right\},
\end{align*}
and let $\mu_{\bH}^1$ be 
the Haar measure on $\bH^1(\A_F)$ given by
\begin{align}\label{muH1DEF}
\int_{\bH^1(\A_F)}f\d \mu_{\bH}^1=
\int_{\bM_{k,d-k}(\A_F)}\int_{G_{d-k}}
f\left(\matr{I_k}Z0{h}\right)
\d \mu^1_{d-k}(h)\d\alpha^{k(d-k)}_F(Z)
\end{align}
for any $f\in\C_c(\bH^1(\A_F))$.
Also set $\Gamma_{\bH}:=\GL_d(F)\cap\bH^1(\A_F)$.
Then
\begin{align}\label{TamagawaHeq1}
\mu_{\bH}^1\bigl(\bH^1(\A_F)/\Gamma_{\bH}\bigr)=1.
\end{align}
Indeed, if 
$\scrF\subseteq \A_F$ is a fundamental domain for $\A_F/F$
and $\scrF_{d-k}\subseteq\GL_{d-k}(\A_F)$
is a fundamental domain for 
$G_{d-k}/\GL_{d-k}(F)$,
then the set
\begin{align*}
\left\{\matr{I_k}{(Z_{i,j})}0h\col Z_{1,1},\, Z_{1,2},\ldots,Z_{k,d-k}\in\scrF\: 
\text{ and } h\in\scrF_{d-k}\right\}
\end{align*}
is a fundamental domain for $\bH^1(\A_F)/\Gamma_{\bH}$,
and so \eqref{TamagawaHeq1}
follows from $\alpha_F(\A_F/F)=1$
and \eqref{TamagawaGLneq1} with $d-k$ in the place of $d$.

Let us fix an arbitrary continuous homomorphism
$t\mapsto\delta_t'=(\delta_{t,\nu}')$,
$\R_{>0}^\times\to\GL_{d-k}(\A_F)$
satisfying $\|\det(\delta_t')\|_{\A_F}=t$ 
and $\delta'_{t,\nu}=I_{d-k}$ for all $t\in\R_{>0}^\times$ and $\nu\in\mathbf{f}$.
(For example we may take $\delta_{t,\nu}'$ diagonal for all $\nu\in\scrP_{F,\infty}$;
we may also let $\delta_{t,\nu}'=I_k$ for all except one
$\nu\in\scrP_{F,\infty}$.)
Then set $\delta_t=\matr{I_k}00{\delta_t'}\in\bH(\A_F)$ for all $t\in\R_{>0}^\times$.
Now we have a
continuous map
$g\mapsto g\,\delta^{-1}_{\|\det(g)\|_{\A_F}}$
from $\bG(\A_F)$ to $G_d$
which induces a homeomorphism between the two quotient spaces
$\bG(\A_F)/\bH(\A_F)$ and $G_d/\bH^1(\A_F)$;
the inverse map is induced by the inclusion
$G_d\hookrightarrow\GL_d(\A_F)$.
Composing with the previous homeomorphism
between $\bG(\A_F)/\bH(\A_F)$ onto $\bU(\A_F)$,
we conclude that the map 
$g\mapsto gp_0$ also induces 
a homeomorphism from $G_d/\bH^1(\A_F)$ onto $\bU(\A_F)$.

Recall that the space $\bM_{d,k}(\A_F)$ is equipped with the measure $\alpha^{dk}_F$;
we write $\alpha^{dk}_F$ also for the restriction of $\alpha_F^{dk}$ to the subset $\bU(\A_F)$.
\begin{lemma}\label{GAFquotLEM2}
The measure $\alpha_F^{dk}$ on $\bU(\A_F)$ corresponds to the quotient measure of $\mu_{\bH}^1$ and $\mu_d^1$
under our fixed homeomorphism between $G_d/\bH^1(\A_F)$ and $\bU(\A_F)$.
\end{lemma}
\begin{proof}
Let $\mu_{\bH}$ be 
the Haar measure on $\bH(\A_F)$ given by
\begin{align}\label{muHDEF}
\int_{\bH(\A_F)}f\d \mu_{\bH}=
\int_{\bM_{k,d-k}(\A_F)}\int_{G_{d-k}}
f\left(\matr{I_k}Z0{h}\right)
\,\d \mu_{d-k}(h) \, \d\alpha_F^{k(d-k)}(Z)
\end{align}
for any $f\in\C_c(\bH(\A_F))$.
Then
for any $f\in\C_c(\bG(\A_F))$,
\begin{align}\label{GAFquotformula1}
\int_{\bU(\A_F)}
\biggl(\int_{\bH(\A_F)}f(\hx h)\d\mu_{\bH}(h)\biggr)\d \alpha_F^{dk}(x)
=\int_{\bG(\A_F)}f(g)\,\|\det(g)\|_{\A_F}^k\d \mu_d(g),
\end{align}
where for each $x\in\bU(\A_F)$,
$\hx$ denotes an arbitrary element in $G_d$ 
satisfying $\hx\, p_0=x\:$ \footnote{Such an $\hx$ exists since
the map $g\mapsto gp_0$ from $G_d$ to $\bU(\A_F)$
is surjective, see our discussion before the statement of the lemma.
Note that for \eqref{GAFquotformula1} to hold we need only require $\hx\in\bG(\A_F)$;
however the requirement $\hx\in G_d$ is important in order to reach \eqref{GAFquotLEM1pf3} below.}
(note that $\int_{\bH(\A_F)}f(\hx h)\d \mu_{\bH}(h)$ only depends on $x$ and not on the choice of $\hx$).
Indeed, in the notation of 
\cite[Ch.\ 2.3]{weil82}, the measure $\mu_d$ as defined in
\eqref{munDEF} equals $\rho_F^{-1}$ times
the Tamagawa measure on $\GL_d(\A_F)$ derived from the
left (and right) invariant gauge-form
$\omega_d=(\det g)^{-d}\d g_{1,1}\d g_{1,2}\,\cdots\, \d g_{d,d}$
on $\GL_d$ by means of the convergence factors
$(\lambda_\nu)$ given by $\lambda_\nu=1$ for $\nu\mid\infty$
and $\lambda_\nu=1-\|\pi_\nu\|_\nu$ for $\nu\in\mathbf{f}$.
Also, in the same notation, $\mu_{\bH}$ defined in \eqref{muHDEF}
equals $\rho_F^{-1}$ times the Tamagawa measure on $\bH(\A_F)$ derived from the
left invariant gauge-form
$\omega_{\bH}=(\det g)^{-(d-k)}\d g_{1,k+1}\d g_{1,k+2}\,\cdots\, \d g_{d,d}$
on $\bH$ by means of the same convergence factors
$(\lambda_\nu)$,
and the measure $\alpha_F^{dk}$ on $\bU(\A_F)$ equals 
the Tamagawa
measure 
derived from the gauge-form
$\omega_{\bU}=\d x_{1,1} \d x_{1,2}\,\cdots\,\d x_{d,k}$
on $\bU$ 
(where $x_{i,j}$ denote the entries of $x\in\bM_{d,k}$)
by means of the convergence factors $(1)$. 
One also verifies that 
$\omega_{\bU}$ is relatively $\bG$-invariant
belonging to the character $g\mapsto \det(g)^k$, 
and 
$\omega_d,\omega_{\bH}$ and $\omega_{\bU}$
match together algebraically.
Hence \eqref{GAFquotformula1} follows from
\cite[Thm.\ 2.4.3]{weil82}.

Now let $t\mapsto\delta_t$ be the one-parameter subgroup of $\bH(\A_F)$
which we fixed in the discussion before the lemma.
We will apply \eqref{omegan1DEF} with this specific choice of $\delta_t$.
Let us also fix
an arbitrary function $f_2\in\C_c(\R^+)$ with $\int_{\R^+}f_2(t)\,t^{k-1}\d t\neq0$.

Now, given any 
$f_1\in \C_c(G_d)$, we define
$f\in\C_c(\bG(\A_F))$ through 
\begin{align*}
f(g)=f_1(g\,\delta_{\|\det(g)\|_{\A_F}}^{-1})\cdot f_2(\|\det(g)\|_{\A_F}).
\end{align*}
Applying \eqref{GAFquotformula1} to this function $f$,
and then using \eqref{muHDEF}
and \eqref{omegan1DEF}, we have
\begin{align}\label{GAFquotLEM1pf2}
\int_{\bU(\A_F)}
\int_{\bM_{k,d-k}(\A_F)}\int_{\GL_{d-k}(\A_F)}
f\left(\hx \matr{I_k}Z0{h}\right)
\d \mu_{d-k}(h)\, d\alpha_F^{k(d-k)}(Z)\d \alpha_F^{dk}(x)
\\\notag
=\int_{G_d} f_1(g)\d \mu_d^1(g)
\cdot \int_{\R^+} f_2(t)\,t^{k-1}\d t.
\end{align}
Furthermore, by \eqref{omegan1DEF} with $d-k$ in place of $d$
and $\delta_t'$ in place of $\delta_t$,
\begin{align*}
\int_{\GL_{d-k}(\A_F)}f\left(\hx \matr{I_k}Z0{h}\right)\d \mu_{d-k}(h)
=\int_{\R^+}\int_{G_{d-k}}f\left(\hx \matr{I_k}Z0{h\delta_t'}\right)\d \mu^1_{d-k}(h)\,\frac{\d t}t,
\end{align*}
and using $\matr{I_k}Z0{h\delta_t'}=\matr{I_k}{Z{\delta_t'}^{-1}}0h \delta_t$,
substituting $Z=Z_{\new}\delta_t'\:$,
and then using \eqref{muH1DEF},
it follows that the left hand side of 
\eqref{GAFquotLEM1pf2} equals
\begin{align}\label{GAFquotLEM1pf3}
\int_{\bU(\A_F)}\int_{\bH^1(\A_F)}f_1(\hx h)\d \mu_{\bH}^1(h)\d \alpha^{dk}_F(x)
\cdot \int_{\R^+} f_2(t)\,t^{k-1}\d t.
\end{align}
Hence, since $\int_{\R^+}f_2(t)\,t^{k-1}\d t\neq0$,
we conclude that
\begin{align}\label{GAFquotLEM1pf1}
\int_{\bU(\A_F)}\int_{\bH^1(\A_F)}f_1(\hx h)\d \mu_{\bH}^1(h)\d \alpha^{dk}_F(x)
=\int_{G_d}f_1(g)\d \mu_d^1(g).
\end{align}
We have proved that this formula holds for all $f_1\in\C_c(G_d)$;
hence \autoref{GAFquotLEM2} is proved.
\end{proof}

Now we apply \cite[Lemma 2.4.2]{weil82}
to the quotient $G_d/\bH^1(\A_F)$
and the discrete subgroups $\Gamma:=\GL_d(F)$ of $G_d$
and $\Gamma_{\bH}:=\Gamma\cap\bH^1(\A_F)$ of $\bH^1(\A_F)$,
making use of \eqref{TamagawaHeq1} and \autoref{GAFquotLEM2}.
It follows that for any 
continuous function
$f:\bU(\A_F)\to\R_{\geq0}$,
\begin{align}
\int_{\bU(\A_F)}f\d \alpha^{dk}_F
=\int_{G_d/\Gamma}\sum_{\gamma\in\Gamma/\Gamma_{\bH}}
f(g\gamma p_0)\d \mu_d^1(g).
\end{align}
Using now the fact that 
$\alpha^{dk}_F(\bM_{d,k}(\A_F)\smallsetminus\bU(\A_F))=0$
by \cite[Lemma 3.4.1]{weil82},
and the fact that when $\gamma$ runs through 
$\Gamma/\Gamma_{\bH}$,
$\gamma p_0$ visits each matrix in $\bM_{d,k}(F)$ with linearly independent columns
exactly once,
we obtain:
\begin{proposition}\label{KIMTHM3p1prop}
For any $1\leq k<d$ and any continuous function $f:\bM_{d,k}(\A_F)\to\R_{\geq0}$,
\begin{align}\label{KIMTHM3p1propRES}
\int_{\scrX_d}\sum_{\substack{X \in \bM_{d,k}(F) \\ \text{columns indep}}} f(gX) \d\mu_d^1(g)
=\int_{\bM_{d,k}(\A_F)}f\d \alpha^{dk}_F.
\end{align}
The relation in \eqref{KIMTHM3p1propRES} is an equality of numbers in $\R_{\geq0}\cup\{+\infty\}$,
i.e., either both sides of \eqref{KIMTHM3p1propRES}
are $+\infty$ or else both are finite and equal.
\end{proposition}

Note that \autoref{KIMTHM3p1prop} gives the same formula
as in the first part of \autoref{thm:arf},
but for a different class of test functions.
Using \autoref{KIMTHM3p1prop} it is now also fairly easy to
prove the second part of \autoref{thm:arf},
for the class of non-negative continuous test functions:
\begin{proposition}\label{adelicRogersPOSITIVEprop}
For any $1\leq k<d$ and any continuous function $f:\bM_{d,k}(\A_F)\to\R_{\geq0}$,
\begin{align}\label{adelicRogersPOSITIVEpropRES}
\int_{\scrX_d}\sum_{\substack{X \in \bM_{d,k}(F) \\ \text{columns nonzero}}} f(gX) \d\mu_d^1(g)
=
\sum_{m=1}^{k}\sum_D\int_{\bM_{d,m}(\A_F)}f(XD)\,\d \alpha^{dm}_F(X),
\end{align}
where the sum over $D$ is over all $m\times k$ row-reduced echelon forms over $F$
of rank $m$ 
with no zero columns.
The relation in \eqref{adelicRogersPOSITIVEpropRES} is again
an equality of numbers in $\R_{\geq0}\cup\{+\infty\}$.
\end{proposition}
(Note that there is exactly one $k\times k$
row-reduced echelon form over $F$ of rank $k$, namely $D=I_k$;
hence the contribution for $m=k$ in the right hand side of
\eqref{adelicRogersPOSITIVEpropRES} is simply
$\int_{\bM_{d,k}(\A_F)}f(X)\,\d \alpha^{dm}_F(X)$.)
\begin{proof}
All equalities in the following proof are equalities between
numbers in $\R_{\geq0}\cup\{+\infty\}$.
By basic linear algebra we have, for any function $f_1:\bM_{d,k}(\A_F)\to\R_{\geq0}$
\begin{align}\label{adelicRogersPOSITIVEproppf1}
\sum_{\substack{X \in \bM_{d,k}(F) \\ \text{columns nonzero}}} f_1(X)
=\sum_{m=1}^k\sum_D
\sum_{\substack{X \in \bM_{d,m}(F) \\ \text{columns indep}}} f_1(XD),
\end{align}
where the sum over $D$ extends over the same set of
$m\times k$ matrices as in \eqref{adelicRogersPOSITIVEpropRES}.
For any $g\in G_d$ we apply \eqref{adelicRogersPOSITIVEproppf1}
to the function $f_1(X):=f(gX)$,
and then integrate the resulting equality over $g$ in a fundamental domain for $G_d/\GL_d(F)$
and change the order between the integration and the summation over $m$ and $D$,
as is permitted by the monotone convergence theorem.
This gives:
\begin{align}
\int_{\scrX_d}\sum_{\substack{X \in \bM_{d,k}(F) \\ \text{columns nonzero}}} f(gX) \d\mu_d^1(g)
=\sum_{m=1}^k\sum_D\int_{\scrX_d}\sum_{\substack{X \in \bM_{d,m}(F) \\ \text{columns indep}}} f(gXD) \d\mu_d^1(g).
\end{align}
Here, for any $m$ and $D$ appearing in the right hand side we have,
by the formula in \autoref{KIMTHM3p1prop}
applied to the continuous function $X\mapsto f(XD)$, $\bM_{d,m}(\A_F)\to\R_{\geq0}$:
\begin{align}\label{adelicRogersPOSITIVEproppf2}
\int_{\scrX_d}\sum_{\substack{X \in \bM_{d,m}(F) \\ \text{columns indep}}} f(gXD) \d\mu_d^1(g)
=\int_{\bM_{d,m}(\A_F)}f(XD)\d \alpha^{dm}_F(X).
\end{align}
Hence we obtain the formula in \eqref{adelicRogersPOSITIVEpropRES}.
\end{proof}

\begin{remark}\label{KIMcomparisonREM1}
The formula in \eqref{adelicRogersPOSITIVEproppf2}
corresponds to \cite[Theorem 3.1]{kim24},
except that there the formula was proved for a different class of test functions $f$.
Note that the class of non-negative continuous functions which we work with here 
is in a certain sense more flexible than the class of test functions 
in \cite[Theorem 3.1]{kim24},
allowing us to deduce the general formula
\eqref{adelicRogersPOSITIVEproppf2} immediately from the special case
appearing in \autoref{KIMTHM3p1prop}.
%
%

The derivation of \autoref{adelicRogersPOSITIVEprop}
from \autoref{KIMTHM3p1prop}
corresponds to the discussion at the end of Section 4.5 in
\cite{kim24};
in particular \eqref{adelicRogersPOSITIVEproppf1} is
a corrected version\footnote{The correction concerns the summation ranges for the sums over $D$.}
of the formula in \cite[p.\ 38 (line 6)]{kim24};
this leads to the correction mentioned in \autoref{thm:arfcorrREM}. 
\end{remark}

In order to prove \autoref{thm:arf} from 
\autoref{KIMTHM3p1prop},
the key step is to prove the following convergence statement.
\begin{proposition}\label{adelicRogersCONVprop}
For any $1\leq k<d$ and any continuous function $f:\bM_{d,k}(\A_F)\to\R_{\geq0}$
of compact support,
\begin{align}\label{adelicRogersCONVpropRES}
\sum_{m=1}^{k}\sum_D\int_{\bM_{d,m}(\A_F)}f(XD)\,\d \alpha^{dm}_F(X)
<\infty.
\end{align}
(The sum over $D$ is as in \autoref{adelicRogersPOSITIVEprop}.)
\end{proposition}
As we have already remarked,
a crucial ingredient in the following proof,
just as in the proof of \cite[Cor.\ 11]{gargavasv23},
is an estimate due to W.\ M.\ Schmidt, \cite[Theorem 3]{schmidt67},
on the number of linear subspaces of $F^k$ of bounded height.
\begin{proof}
Let $n=[F:\Q]$, and fix a $\Z$-basis $f_1,\ldots,f_n$ of $\scrO_F$;
this is then also a 
basis of $F$ as a vector space over $\Q$.
For any $x\in F$ we write $[x]\in\Q^n$ for the coordinates of
$x$ in 
this basis,
i.e., $[x]=(c_1,\ldots,c_n)\in\Q^n$ when $x=\sum_{i=1}^n c_if_i$.
Let the map $x\mapsto\tx$, $F\to\bM_n(\Q)$, 
be the corresponding regular representation,
so that $[xy]=[x]\,\ty$ for all $x,y\in F$.


For each place $p\in\scrP_{\Q}$ there is a standard isomorphism
of $\Q_p$-algebras $F\otimes_{\Q}\Q_p\cong\prod_{\nu\mid p}F_{\nu}$,  
and using the basis $f_1,\ldots,f_n$ we also have an isomorphism of
$\Q_p$-modules $\Q_p^n\cong F\otimes_{\Q}\Q_p$,
\mbox{$(c_1,\ldots,c_n)\mapsto\sum_{i=1}^n f_i\otimes c_i$.}
Using these isomorphisms for each quasifactor in 
$\A_F= 
\prod_{p\in\scrP_{\Q}}'\bigl(\prod_{\nu\mid p}F_\nu\bigr)$
we obtain an isomorphism of topological $\A_{\Q}$-modules 
$J:\A_F\stackrel{\sim}{\rightarrow}\A_{\Q}^n$
\cite[Theorem IV-1.1]{weil95}.
It follows from the construction that 
\begin{align}\label{adelicRogersCONVproppf2}
J(xy)=J(x)\,\ty,\qquad\forall x\in \A_F,\: y\in F.
\end{align}
Furthermore, the Haar measure $\alpha_F$ corresponds to a multiple
of $\alpha_{\Q}^d$ under $J$. 

By abuse of notation, for any $d,m$ let us write $J$ also for the
isomorphism of topological $\A_{\Q}$-modules
$J:\bM_{d,m}(\A_F)\stackrel{\sim}{\rightarrow}\bM_{d,mn}(\A_{\Q})$
which maps $X=(x_{i,j})_{\substack{i=1,\ldots,d\\ j=1,\ldots,m}}$
in $\bM_{d,m}(\A_F)$ to the block matrix
$(J(x_{i,j}))_{\substack{i=1,\ldots,d\\ j=1,\ldots,m}}$
with each $J(x_{i,j})\in\A_{\Q}^n$ viewed as a $1\times n$ matrix.
Also 
for any $D=(d_{i,j})$ in $\bM_{m,k}(F)$ we write
$\tD:=(\td_{i,j})\in\bM_{mn,kn}(\Q)$.
It then follows from \eqref{adelicRogersCONVproppf2} that 
\begin{align}
J(XD)=J(X)\tD
\qquad \forall X\in\bM_{d,m}(\A_F),
\: D\in\bM_{m,k}(F).
\end{align}

Using the isomorphism
$J:\bM_{d,k}(\A_F)\stackrel{\sim}{\rightarrow}\bM_{d,kn}(\A_{\Q})$,
it follows that it suffices to prove that
for any continuous function $f:\bM_{d,kn}(\A_\Q)\to\R_{\geq0}$
of compact support,
\begin{align}\label{adelicRogersCONVproppf3}
\sum_{m=1}^{k}\sum_D
\int_{\bM_{d,mn}(\A_\Q)}f(X\tD)\,\d \alpha^{dmn}_{\Q}(X)
<\infty.
\end{align}
Given such a function $f$, 
there exists a compact set $C\subseteq\A_\Q$
such that 
\begin{align*}
\supp(f)\subseteq C^{dkn}:=\{X=(X_{i,j})\in\bM_{d,kn}(\A_\Q)\col X_{1,1},\: X_{1,2},\: \cdots, X_{d,kn}\in C\}.
\end{align*}
Since $f$ is bounded, it now suffices to prove
that \eqref{adelicRogersCONVproppf3} holds with 
the characteristic function of $C^{dkn}$ in place of $f$,
viz.,
\begin{align}\label{adelicRogersCONVproppf4}
\sum_{m=1}^{k}\sum_D
\alpha^{dmn}_{\Q}\bigl(\bigl\{
X\in\bM_{d,mn}(\A_\Q)
\col X\tD\in C^{dkn}\bigr\}\bigr)
<\infty.
\end{align}
For each place $p$ of $\Q$, let $C_p$ be the image of $C$ under
the projection map from $\A_\Q$ to $\Q_p$;
then each $C_p$ is compact and $C_p\subseteq\Z_p$ for almost all finite places $p$.
Also $C\subseteq\prod_p C_p$, and 
hence it suffices to prove that 
\eqref{adelicRogersCONVproppf4} holds with $\prod_p C_p$
in place of $C$;
thus we will from now on assume that $C=\prod_p C_p$.

Fix some $m$ and $D$ appearing in the sum in 
\eqref{adelicRogersCONVproppf4}.
Then
\begin{align}\label{adelicRogersCONVproppf8}
\alpha^{dmn}_{\Q}\bigl(\bigl\{
X\in\bM_{d,mn}(\A_\Q)
\col X\tD \in C^{dkn}\bigr\}\bigr)
=\alpha^{mn}_{\Q}\bigl(\bigl\{\vecx\in\A_\Q^{mn}\col \vecx\,\tD \in C^{kn}\bigr\}\bigr)^d.
\end{align}
For any $T\in\GL_{mn}(\Q)$,
the map $\vecx\mapsto \vecx T$ preserves
$\alpha_{\Q}^{mn}$
since $\|\det(T)\|_{\A_{\Q}}=1$;
hence
\begin{align}\label{adelicRogersCONVproppf6}
\alpha^{mn}_{\Q}\bigl(\bigl\{\vecx\in\A_\Q^{mn}\col \vecx\,\tD \in C^{kn}\bigr\}\bigr)
&=\alpha^{mn}_{\Q}\bigl(\bigl\{\vecx\in\A_\Q^{mn}\col \vecx\,T\,\tD \in C^{kn}\bigr\}\bigr)
\\\notag
&=\prod_{p\in\scrP_\Q}\alpha_p^{mn}\bigl(\bigl\{\vecx\in\Q_p^{mn}
\col \vecx\, T\,\tD \in C_p^{kn}\bigr\}\bigr),
\end{align}
where the last equality holds since 
$C=\prod_p C_p$ and using \eqref{alphaFDEF} and $\Delta_\Q=1$.
Set $V_D:=\Q^{mn}\tD $;
this is a $\Q$-linear subspace of $\Q^{kn}$
of dimension $mn$;
indeed, it is the image of the $m$-dimensional $F$-linear subspace
$F^mD$ of $F^k$ under the $\Q$-linear bijection
$(x_1,\ldots,x_k)\mapsto([x_1],\ldots,[x_k])$ 
from $F^k$ onto $\Q^{kn}$.
Hence $\Lambda_D:=V_D\cap\Z^{kn}$
is a free $\Z$-module of rank $mn$.
Let us fix a $\Z$-basis $\vecb_1,\ldots,\vecb_{mn}$ of $\Lambda_D$,
and choose $T\in\GL_{mn}(\Q)$ so that
the matrix 
$T\,\tD $ 
has row vectors
$\vecb_1,\ldots,\vecb_{mn}$,    
in this order.
Then for each finite place $p$ we have
\begin{align}\label{adelicRogersCONVproppf5}
\bigl\{\vecx\in\Q_p^{mn}\col \vecx\,T\,\tD \in \Z_p^{kn}\bigr\}
=\bigl\{(x_1,\ldots,x_{mn})\in\Q_p^{mn}\col {\textstyle\sum_{i=1}^{mn}x_i\vecb_i}\in \Z_p^{kn}\bigr\}
=\Z_p^{mn}.
\end{align}
Here the second equality is proved as follows.
It is immediate that $\sum_{i=1}^{mn}x_i\vecb_i\in\Z_p^{kn}$
holds for every $\vecx\in\Z_p^{mn}$.
Conversely, assume that
$\vecx\in\Q_p^{mn}$ satisfies
$\sum_{i=1}^{mn}x_i\vecb_i\in\Z_p^{kn}$.
Take $\ell\in\Z_{\geq0}$,
$\vecx'\in\Z^{mn}$ and $\vecx''\in\Z_p^{mn}$ 
so that $\vecx=p^{-\ell}\vecx'+\vecx''$.
Then 
$\sum_{i=1}^{mn}x_i''\vecb_i\in\Z_p^{kn}$,
and hence
$p^{-\ell}\sum_{i=1}^{mn}x_i'\vecb_i\in\Z_p^{kn}$.
But the vector
$p^{-\ell}\sum_{i=1}^{mn}x_i'\vecb_i$
also belongs to $p^{-\ell}\Z^{kn}$ and to $V_D$,
and since $p^{-\ell}\Z\cap\Z_p=\Z$,
it follows that
$p^{-\ell}\sum_{i=1}^{mn}x_i'\vecb_i\in V_D\cap\Z^{kn}=\Lambda_D$.
But $\vecb_1,\ldots,\vecb_{mn}$ is a $\Z$-basis of $\Lambda_D$;
thus $p^{-\ell}x'_i\in\Z$ for all $i$, 
viz., $p^{-\ell}\vecx'\in\Z^{mn}$,
and hence $\vecx\in\Z_p^{mn}$.
This completes the proof of 
the second equality in
\eqref{adelicRogersCONVproppf5}.

For each finite place $p$, since $C_p$ is compact,
$C_p$ can be covered by a finite number $N_p$ of translates of $\Z_p$,
and then $C_p^{kn}$ can be covered by $N_p^{kn}$
translates of $\Z_p^{kn}$,
i.e., sets of the form $\veca+\Z_p^{kn}$ with $\veca\in\Q^{kn}$.
For each such translate,
by linearity and by \eqref{adelicRogersCONVproppf5},
the set $\{\vecx\in\Q_p^{mn}\col \vecx\,T\,\tD \in \veca+\Z_p^{kn}\bigr\}$
is either empty or a translate of $\Z_p^{mn}$.
Hence
\begin{align}
\alpha_p^{mn}\bigl(\bigl\{\vecx\in\Q_p^{mn}
\col \vecx\, T\,\tD \in C_p^{kn}\bigr\}\bigr)
\leq N_p^{kn}.
\end{align}
Hence by \eqref{adelicRogersCONVproppf6}, 
since we can take $N_p=1$ for almost all $p$,
and since $\alpha_\infty^{mn}=\vol$,
the standard Lebesgue volume measure on $\R^{mn}$,
we have
\begin{align}\label{adelicRogersCONVproppf7}
\alpha^{mn}_{\Q}\bigl(\bigl\{\vecx\in\A_\Q^{mn}\col \vecx\,\tD \in C^{kn}\bigr\}\bigr)
\ll \vol\bigl(\bigl\{\vecx\in\R^{mn}\col 
{\textstyle\sum_{i=1}^{mn}x_i\vecb_i}\in C_\infty^{kn}\bigr\}\bigr),
\end{align}
where the implied constant depends only on the compact set $C=\prod_p C_p$.
Now take $R>0$ so that $C_\infty^{kn}$ is contained in the ball of radius $R$
in $\R^{kn}$ centered at the origin.
Then the volume in \eqref{adelicRogersCONVproppf7} is
bounded above by $V/\mathsf{d}(\Lambda_D)$
where $V$ is the volume of an $mn$-dimensional ball of radius $R$
and $\mathsf{d}(\Lambda_D)$ is the determinant of $\Lambda_D$
viewed as an $mn$-dimensional lattice in $\R^{kn}$,
i.e., the $mn$-dimensional volume of the 
parallelotope spanned by $\vecb_1,\ldots,\vecb_{mn}$
(notation as in W.\ M.\ Schmidt, \cite[\S 2]{schmidt67}).
Hence, recalling \eqref{adelicRogersCONVproppf4} and \eqref{adelicRogersCONVproppf8},
it now suffices to prove that
\begin{align}\label{adelicRogersCONVproppf9}
\sum_{m=1}^{k}\sum_D
\mathsf{d}(\Lambda_D)^{-d}
<\infty.
\end{align}

Next let us note that for any $m$ and $D$ as in the above sum,
$\Lambda_D$ viewed as a lattice in $\R^{kn}$
equals the image of 
$F^mD\cap\scrO_F^k$
under the map $(x_1,\ldots,x_k)\mapsto([x_1],\ldots,[x_k])$ 
from $F^k$ to $\R^{kn}$.
Hence we have
\begin{align}\label{adelicRogersCONVproppf10}
\mathsf{d}(\Lambda_D)\asymp H'(F^mD),
\end{align}
where $H'(F^mD)$ is the \textit{height} of the linear subspace $F^mD$ of $F^k$,
as defined in \cite[eq.\ (8)]{schmidt67},
and with the implied constants in \eqref{adelicRogersCONVproppf10}
only depending on $F$ and $k$.\footnote{(Indeed, let $\beta_1$ be the $\R$-linear map from $\R^n$ to $\R^n$
defined by $\beta_1(\vece_i)=(f_i^{[1]},\ldots,f_i^{[n]})$ for $i=1,\ldots,n$,\label{footnoteonasympHpFmD}
where we are using the notation in 
\cite[eq.\ (7 bis)]{schmidt67}, 
and where $\vece_i=(0,\ldots,1,\ldots,0)$, the $i$th standard basis vector.
Define $\beta:\R^{kn}\to\R^{kn}$ by $\beta(\vecv_1,\ldots,\vecv_k)=(\beta_1(\vecv_1),\ldots,\beta_1(\vecv_k))$.
Then the composition of our map $F^k\to\R^{kn}$
and the map $\beta$ equals W.\ M.\ Schmidt's map $\rho$ \cite[p.\ 435(mid)]{schmidt67},
and now our claim follows from the fact that,
since $\beta$ is non-singular,
the volume scaling factor of $\beta$ restricted to any linear subspace $U$ of
$\R^{kn}$ is bounded from above and below by positive constants independent of $U$.)}
Hence it now suffices to prove that for each fixed $m\in\{1,\ldots,k\}$,
\begin{align}\label{adelicRogersCONVproppf11}
\sum_U H'(U)^{-d}<\infty,
\end{align}
where the sum is taken over all $m$-dimensional $F$-linear subspaces $U$ of $F^k$.
This is obvious if $m=k$ (since then the sum has only one term); 
hence from now on we assume $m<k$.

Let us write $N(x)$ for the number of 
$m$-dimensional $F$-linear subspaces $U\subseteq F^k$
satisfying $H'(U)<x$.
Then by \cite[Theorems 1 and 3]{schmidt67},
we have $N(x)\ll x^k$ for all $x\geq1$;
also by \cite[Lemma 8]{schmidt67},
$N(1)=0$.
Now by dyadic decomposition, the sum in 
\eqref{adelicRogersCONVproppf11} is bounded above by
\begin{align*}
\sum_{j=0}^\infty \bigl(N(2^{j+1})-N(2^{j})\bigr)\cdot 2^{-jd}
\leq \sum_{j=0}^\infty N(2^{j+1})\cdot 2^{-jd}
\ll\sum_{j=0}^\infty 2^{(j+1)k-jd},
\end{align*}
and the last sum is finite since $k<d$.
This completes the proof of \autoref{adelicRogersCONVprop}.
\end{proof}

\begin{proof}[Proof of \autoref{thm:arf}]
This is a standard deduction from \autoref{adelicRogersPOSITIVEprop}
and \autoref{adelicRogersCONVprop},
using basic integration theory.
Indeed, write $Y=\bM_{d,k}(\A_F)$,
define the Borel measures 
$\lambda_1$ and $\lambda_2$
on $Y$
by setting, for any Borel set $A\subseteq Y$:
\[
    \lambda_1(A)=\int_{\scrX_d}\sum_{\substack{X \in \bM_{d,k}(F) \\ \text{columns nonzero}}} \chi_A(gX) \d\mu_d^1(g),
\]
and
\[
    \lambda_2(A)=\sum_{m=1}^{k}\sum_D\int_{\bM_{d,m}(\A_F)}\chi_A(XD)\,\d \alpha^{dm}_F(X).
\]
Then for any Borel measurable function
$f:Y\to\R_{\geq0}$,
by using the definition of $\int_Y f\d\lambda_j$ and the Monotone Convergence Theorem,
we have
\begin{align}\notag 
\int_{Y} f\d\lambda_1
=\int_{\scrX_d}\sum_{\substack{X \in \bM_{d,k}(F) \\ \text{columns nonzero}}} f(gX) \d\mu_d^1(g)
\end{align}
and 
\begin{align}\notag 
\int_{Y} f \d\lambda_2
=\sum_{m=1}^{k}\sum_D\int_{\bM_{d,m}(\A_F)}f(XD)\,\d \alpha^{dm}_F(X).
\end{align}

Now by \autoref{adelicRogersPOSITIVEprop}
and \autoref{adelicRogersCONVprop},
we have $\lambda_1(f)=\lambda_2(f)<\infty$ for all
non-negative functions $f\in \C_c(Y)$.
For any compact set $K\subseteq Y$, there exists a non-negative
$f\in \C_c(Y)$
satisfying $\chi_K\leq f$, 
and thus $\lambda_j(K)\leq\int_Y f\d\lambda_j<\infty$
for $j=1,2$.
Hence by \cite[Theorem 2.18]{rudin87},
both $\lambda_1$ and $\lambda_2$ are regular.
Now by the uniqueness argument in the proof of 
the Riesz Representation Theorem
(see, e.g., \cite[p.\ 41]{rudin87}),
it follows that $\lambda_1=\lambda_2$.
Hence \autoref{thm:arf} holds.
\end{proof}

\newpage

\section{Diophantine approximation with congruence conditions}\label{sec:dc}
In this section we will prove \autoref{thm:cdioc}.
We begin with a brief discussion of the methods of \cite{alamghoshyu21} in relation to their main result \autoref{thm:dioccagy}.
Recall from the introduction \S\ref{ssec:cdioc} that \cite{alamghoshyu21} translates diophantine approximation 
over $\mathbb{Q}$ with congruence conditions to a lattice point counting problem with $\mathcal{Y}_{N, \bm{v}}$ being the suitable space of lattices.
The main ingredients of their proof are as follows:
\begin{enumerate}[label=(\Alph*),font=\normalfont,before=\normalfont]
    \item The first moment and a variance bound of the lattice point counting function (mentioned below respectively):
        \label{inga}
        \begin{theorem}[\cite{MarklofStrombergsson2010}, Proposition 7.1]\label{thm:affinesiegel}
            Let $f \in L^1(\mathbb{R}^{d})$.
            Then
            \[
                \int_{\mathcal{Y}_{N, \bm{v}}} \sum_{\bm{w} \in {\left(\mathbb{Z}^{d} + \frac{\bm{v}}{N}\right)} \smallsetminus \{\bm{0}\}} f(g\bm{w}) \d{\mu(g)} = \int_{\mathbb{R}^{d}} f(\bm{x}) \d{\bm{x}}.
            \]
        \end{theorem}
        \begin{theorem}[\cite{ghoshkelmeryu2020}, Corollary 3.4]\label{thm:affinevarbdd}
            For bounded Borel measurable $A \subseteq \mathbb{R}^{d} \smallsetminus \{\bm{0}\}$
            \[
                \int_{\mathcal{Y}_{N, \bm{v}}} {\left|\#\left(g{\left(\mathbb{Z}^{d} + \frac{\bm{v}}{N}\right)} \cap A\right) - \vol(A)\right|}^2 \d{\mu(g)} \ll_{d, N} \vol(A).
            \]
        \end{theorem}
    \item A technique due to W.\ M.\ Schmidt \cite{schmidt60b} which utilizes the above two theorems to prove:
        \label{ingb}
        \begin{theorem}\label{thm:dioccgen}
            Let $\{A_T\}_{T > 1}$ be an increasing family of bounded Borel subsets of $\mathbb{R}^{d} \smallsetminus \{\bm{0}\}$ such that $\lim\limits_{T \to \infty} \vol(A_T) = \infty$.
            Then for all $\varepsilon > 0$ and for $\mu$-a.e.\ $g \in \SL_d( \mathbb{R})$
            \[
                \#{\left(g{\left(\mathbb{Z}^{d} + \frac{\bm{v}}{N}\right)} \cap A_T\right)} = \vol(A_T) + O_\varepsilon{\left(\vol(A_T)^{1/2 + \varepsilon}\right)}.
            \]
        \end{theorem}
    \item An approximation argument \cite[\S3]{alamghoshyu21} to prove \autoref{thm:dioccagy} from \autoref{thm:dioccgen}.
        (We remark that the error term of \autoref{thm:dioccgen} is lost during this approximation.)
        \label{ingc}
\end{enumerate}

Our proof of \autoref{thm:cdioc} will contain analogues of
each of the above steps (A)--(C), but in our different set-up
where, in particular, $\mathcal{Y}_{N, \bm{v}}$
is replaced by the homogeneous adele space $\scrX_d=G_d/\GL_d(F)$.

\subsection{Diophantine approximation over $\CC$ with congruence conditions}\label{ssec:dcx}
We now start the proof of \autoref{thm:cdioc}.
Our first step will be to   
translate the problem of diophantine approximation on $\mathbb{C}$ with congruence conditions 
into a lattice point counting problem over the adeles.
Let $m,n,d,F,\scrO,\frakI,\vecv$ and $\psi$ be as in the statement of
\autoref{thm:cdioc},
and recall that $N(\vartheta,T,\vecv, \mathfrak{I})$ denotes the number of
pairs $(\bm{p}, \bm{q})\in \mathcal{O}^{m} \times \mathcal{O}^{n}$ 
satisfying the inequalities and congruence condition in \eqref{eq:cdioc}.


Let us define the function 
$u : \M_{m \times n}(\mathbb{C}) \to G_d$
through
\begin{align*}
    u(\vartheta)_\nu := \begin{cases}
        \begin{pmatrix}
            1_m & \vartheta \\
            0   & 1_n
        \end{pmatrix}, \ \text{for} \ \nu = \infty \\
        1_d, \ \text{for all} \ \nu \ \text{non-archimedean}.
    \end{cases}
\end{align*}
The for each $\vartheta \in \M_{m \times n}(\mathbb{C})$ 
we define
\begin{align*}
\Lambda_\vartheta:=u(\vartheta)F^{d} \in \mathcal{X}_d.
\end{align*}

Recall that there is a bijective correspondence
$\nu \leftrightsquigarrow \mathfrak{p}$
between the finite places of $F$
and the (non-zero) prime ideals of $\scrO$:
For each prime ideal $\frakp$ of $\scrO$,
we let $\|\cdot\|_{\frakp}$ be the unique absolute value on 
$F$
satisfying $\|\pi\|_{\frakp}=N(\mathfrak{p})^{-1}$ for all $\pi\in\frakp\smallsetminus\frakp^2$,
where $N(\frakp)=\#(\scrO/\frakp)$ is the norm of $\frakp$;
then $\|\cdot\|_{\frakp}=\|\cdot\|_\nu$
for a unique finite place $\nu\in\vecf$
(see \S\ref{notationSEC}).
Thus, sometimes we will write $\norm{ \cdot }_\nu$ and $\norm{ \cdot }_\mathfrak{p}$ interchangeably, and denote the completion of $F$ with respect to this valuation by $F_\nu$ or $F_\mathfrak{p}$.
The closure of $\mathcal{O}$ in this completion will be denoted as $\mathcal{O}_\nu$ or $\mathcal{O}_\mathfrak{p}$ (thus, for us, $\mathcal{O}_\mathfrak{p}$ is not the localization of $\mathcal{O}$ at $\mathfrak{p}$).
We will think of the archimedean valuation on $F$ as corresponding to a `prime' denoted by $\infty$.
An element of $\mathbb{A}_F^{k}$ can thus be thought of as a tuple indexed by the primes of $\mathcal{O}$ 
together with the symbol $\infty$.
We will write `$\mathfrak{p} < \infty$' to denote that $\frakp$ is a genuine prime ideal of $\scrO$.
For any prime $\frakp$ (including $\frakp=\infty$),
by a mild abuse of notation we will write $\|\cdot\|_{\frakp}$ also
for the induced supremum norm on $F_{\frakp}^k$ for any $k\in\Z^+$;
furthermore, for any $\vecx\in\A_F^k$ we will use
``$\|\vecx\|_{\frakp}$'' as a short-hand for $\|\vecx_{\frakp}\|_{\frakp}$.

Let the unique prime factorization of $\frakI$ be 
\begin{align*}
\mathfrak{I} = \prod_{\frakp<\infty} \mathfrak{p}^{\alpha_{\mathfrak{p}}}.
\end{align*}
(Thus $\alpha_{\frakp}\in\Z_{\geq0}$ for all $\frakp$
and $\alpha_{\frakp}=0$ for almost every $\frakp$.)
Next, for any $T>1$, we define $E_T$ to be the following subset of
$\A_F^m\times\A_F^n=\A_F^d$:
\[
    E_{T} 
:= \left\{\bm{z} = (\bm{x}, \bm{y}) \in \mathbb{A}_F^m \times \mathbb{A}_F^{n} : \substack{\displaystyle \norm{\bm{x}}_\infty^m \leq \psi(\norm{\bm{y}}_{\infty}^{n}),\: 1 \leq \norm{\bm{y}}_{\infty}^{n} < T \: \text{ and} \\[5pt] \displaystyle \norm{\bm{z} - \bm{v}}_\mathfrak{p} \leq N(\mathfrak{p})^{- \alpha_\mathfrak{p}} \: \text{ for all } \: \mathfrak{p} < \infty}\right\}.
\]
Since $\bm{v} \in \mathcal{O}^{d}$, we have the following equivalence:
\begin{equation}\label{eq:ftoco}
    \bm{z} \in F^{d} \ \text{and} \ \norm{\bm{z} - \bm{v}}_\mathfrak{p} \leq N(\mathfrak{p})^{- \alpha_\mathfrak{p}} \ \text{for all} \ \mathfrak{p} < \infty \quad\iff\quad \bm{z} \in \mathcal{O}^{d} \ \text{and} \ \bm{z} = \bm{v}\:\: (\text{mod } {\mathfrak{I}}).
\end{equation}
Therefore, for any $\vartheta\in\M_{m \times n}(\mathbb{C})$,
the number of solutions to \eqref{eq:cdioc} can be expressed as
\begin{align}\label{NtTvascounting}
\mathcal{N}(\vartheta, T, \bm{v}, \mathfrak{I})=
\#(E_T \cap \Lambda_\vartheta).
\end{align}


Note that
\begin{align*}
E_T = E_{T, \infty} \times \prod_{\mathfrak{p} < \infty} (\bm{v} + \mathfrak{p}^{ \alpha_\mathfrak{p}} \mathcal{O}_\mathfrak{p}^{d}),
\end{align*}
where
\[
    E_{T, \infty} := \{(\bm{x}, \bm{y}) \in \mathbb{C}^{m} \times \mathbb{C}^{n} : \norm{\bm{x}}_\infty^m \leq \psi(\norm{\bm{y}}_\infty^n) \ \text{and} \ 1 \leq \norm{\bm{y}}_\infty^n < T\}.
\]
Hence by \eqref{alphaFDEF},
\begin{align}\label{alphaFdETformula}
    \alpha_F^d(E_T) = |\Delta_F|^{-d/2}\,\alpha_\infty^d(E_{T, \infty}) \prod_{\mathfrak{p} < \infty} N(\mathfrak{p})^{- \alpha_\mathfrak{p}d} 
= 2^d |\Delta_F|^{-d/2} N(\mathfrak{I})^{-d}\vol(E_{T, \infty}).
\end{align}
(Indeed, recall that ``$\vol$'' denotes the standard $2d$-dimensional volume measure on $\mathbb{C}^d=(\R\oplus i\R)^d$;
thus $\vol=2^{-d}\alpha_\infty^d$; see \S\ref{notationSEC}.)

Furthermore,
\begin{align}\notag
    \vol(E_{T, \infty}) &= \int_{\{\bm{y} \in \mathbb{C}^{n}: 1 \leq \norm{\bm{y}}^{n}_{\infty} < T\}} \int_{\{\bm{x} \in \mathbb{C}^{m} : \norm{\bm{x}}^{m}_{\infty} < \psi({\norm{\bm{y}}}^{n}_{\infty})\}} \d{\bm{x}} \d{\bm{y}} 
\hspace{100pt}
\\\label{eq:etinfvolPRE}
                        &= \pi^m \int_{\{\bm{y} \in \mathbb{C}^{n} : 1 \leq {\norm{\bm{y}}}^{n}_{\infty} < T\}} \psi({\norm{\bm{y}}}^{n}_{\infty}) \d{\bm{y}}.
\end{align}
To compute the last integral,
let us extend $\psi$ by setting $\psi(r):=0$ for all $0\leq r<1$.
Substituting $\vecy=\bigl(\sqrt{t_1}e^{i\omega_1},\ldots,\sqrt{t_n} e^{i\omega_n}\bigr)$,
and making use of the symmetry under permutation of the $n$ coordinates,
we see that the last integral in \eqref{eq:etinfvolPRE} equals
\begin{align*}
\pi^n\, n!\, \int_{0<t_1<\cdots<t_n<T} \psi(t_n^n)\,\d t_1\,\cdots\,\d t_n
=\pi^n\, n!\,\int_0^T\psi(t_n^n)\,\frac{t_n^{n-1}}{(n-1)!}\, \d t_n
=\pi^n\int_0^T\psi(t)\, \d t.
\end{align*}
Hence we conclude:
\begin{equation}\label{eq:etinfvol}
    \vol(E_{T, \infty}) = \pi^{d} \int_{1}^{T} \psi(t) \d{t}.
\end{equation}

In view of \eqref{NtTvascounting}, \eqref{alphaFdETformula} 
and \eqref{eq:etinfvol},
\autoref{thm:cdioc} can be restated as follows:
\begin{theorem}\label{thm:cdiocnew}
Assume that $\int_{1}^{\infty} \psi(t) \d{t} = \infty$.
Then for almost every $\vartheta \in \M_{m \times n}(\mathbb{C})$,
    \[
        \#(E_T \cap \Lambda_\vartheta) 
\thicksim \alpha_F^d(E_T), \qquad \text{as } \: T \to \infty.
    \]
\end{theorem}

\subsection{Analogue of Schmidt \cite{schmidt60a}}\label{ssec:schmidt60aanalogue}
In this section we state an analogue of \autoref{thm:dioccgen}.
The following variance bound is a crucial input in this theorem.

\begin{theorem}[Variance bound; Theorem 1.3 in \cite{kim24}]\footnotemark
\label{thm:vbac}
Let $d\geq3$.
For each place $\nu$ of $F$, let $A_\nu$ be a Borel measurable subset of $F_\nu^d$,
with $A_\nu\subseteq\scrO_\nu^d$ for all except finitely many $\nu<\infty$, 
and set $A=\prod_{\nu\in\scrP_F}A_\nu\subseteq\A_F^d$.
Then 
\begin{align}\label{thm:vbacRES}
\int_{\mathcal{X}_d} \bigl(\#(A\cap \Lambda\smallsetminus\{\bn\})-\alpha_F^d(A)\bigr)^2\,\d\mu_d^1(\Lambda) \ll \alpha_F^d(A),
\end{align}
where the implied constant depends only on $F$ and $d$. 
\end{theorem}

\footnotetext{We apply
\cite[Theorem 1.3]{kim24} with ``$f_\infty$'' being the characteristic function of the set $A_\infty$.
Note that since $F$ is an imaginary quadratic number field,
the condition 
(1.2) 
in the statement of \cite[Theorem 1.3]{kim24}
is automatically fulfilled, with $C=1$.
It should also be noted that 
$\int_{\mathcal{X}_d} \#(A\cap \Lambda\smallsetminus\{\bn\})\,\d\mu_d^1(\Lambda)=\alpha_F^d(A)$,
by \autoref{thm:arf} applied with $k=1$ and $f=$ the characteristic function of $A$.}

We stress that \cite[Theorem 1.3]{kim24} is more general than
\autoref{thm:vbac} in that it applies for a general number field,
and admits more general integrands.
A crucial ingredient in the proof in \cite{kim24} 
is the adelic Rogers' formula, \autoref{thm:arf}
= \cite[Theorem 1.2]{kim24},
applied with $k=1$ and $k=2$.
This is also the reason why we have to assume $d\geq3$
in the statement of \autoref{thm:vbac}.   

Let us now apply
\autoref{thm:vbac}
with $A_{\frakp}=\bm{v} + \mathfrak{p}^{ \alpha_\mathfrak{p}}\mathcal{O}_\mathfrak{p}^{d}$
for each finite place $\frakp$,
and with $A_\infty$ being a subset of 
$\CC^d\smallsetminus\{\bn\}$.
Then in \eqref{thm:vbacRES} we have 
$A\cap \Lambda\smallsetminus\{\bn\}=A\cap \Lambda$
for all $\Lambda\in\scrX_d$,
and also, 
as in \eqref{alphaFdETformula},
\begin{align}\label{alphaFdAformula}
    \alpha_F^d(A) = 2^d |\Delta_F|^{-d/2} N(\mathfrak{I})^{-d}\vol(A_{\infty}).
\end{align}
With these facts, Schmidt's argument in \cite[Theorem 1]{schmidt60a} can be applied
to deduce the following:
\begin{theorem}\label{thm:diocca}
    Let $\{A_{T, \infty}\}_{T > 1}$ be an increasing family of bounded Borel subsets of $\mathbb{C}^{d} \smallsetminus \{\bm{0}\}$ such that $\lim\limits_{T \to \infty} \vol(A_{T, \infty}) = \infty$.
    In addition, let $A_T = A_{T, \infty} \times \prod\limits_{\mathfrak{p} < \infty} (\bm{v} + \mathfrak{p}^{ \alpha_\mathfrak{p}} \mathcal{O}_\mathfrak{p}^{d})$ with $\alpha_{\mathfrak{p}}$ and $\mathfrak{I}$ as above.
    Then for all $\varepsilon > 0$ and for $\mu_d^1$-almost every $g \in G_d$,
    \begin{equation}\label{eq:schasymp}
        \#\left(A_T \cap gF^{d}\right) = \alpha_F^d(A_T) + O\Bigl(\alpha_F^d(A_T)^{1/2 + \varepsilon}\Bigr)
\qquad\text{as }\: T\to\infty.
    \end{equation}
\end{theorem}
\begin{proof}
We here prove \autoref{thm:diocca} as an application of
\cite[Lemma 10]{vS79},
which is abstracted from the work of Schmidt.
Without loss of generality, we may assume that
the family $\{A_{T,\infty}\}$ is parametrized by
$T\in[1,\infty)$, and that we have
$\alpha_F^d(A_T)=T$ whenever $T$ is an integer.\footnote{This claim is verified by
first, if necessary, throwing away $A_{T,\infty}$ for $T$ small;
and then reparametrizing and possibly enlarging the family
$\{A_{T,\infty}\}$, making use of 
\eqref{alphaFdAformula} and \cite[Lemma 1]{schmidt60a}.}
We will apply \cite[Lemma 10]{vS79}
to the probability space space $(\scrX_d,\mu_d^1)$
and with the functions $f_1,f_2,\ldots$ defined through
$f_T(\Lambda):=\#((A_T\smallsetminus A_{T-1})\cap\Lambda)$ for $T\in\Z^+$,
where we set $A_0:=\emptyset$.
Note that for any integers $0\leq m<n$, we have
\begin{align*}
\int_{\scrX_d} &\biggl(\sum_{m<k\leq n}f_k(\Lambda)-(n-m)\biggr)^2\,d\mu_d^1(\Lambda)
\\
&=\int_{\scrX_d}\bigl(\#((A_n\smallsetminus A_m)\cap\Lambda)-\alpha_F^d(A_n\smallsetminus A_m)\bigr)^2\,d\mu_d^1(\Lambda)
\ll n-m,
\end{align*}
by \autoref{thm:vbac}
(which applies since $A_n\smallsetminus A_m=(A_{n,\infty}\smallsetminus A_{m,\infty})\times\prod_{\frakp<\infty}A_{\frakp}$).
This shows that the assumption in
\cite[Lemma 10]{vS79}
is fulfilled\footnote{With all the \textit{numbers} ``$f_k$'' and ``$\varphi_k$''
in the statement of  \cite[Lemma 10]{vS79} taken to equal $1$.}.
The conclusion from that lemma is that,
for $\mu_d^1$-almost every $\Lambda\in\scrX_d$, we have
\begin{align}\label{thm:dioccapf1}
\#(A_n\cap\Lambda)
=\sum_{k=1}^nf_k(\Lambda)=n+O(n^{1/2+\ve})\qquad
\text{as }\: n\to\infty\text{ through $\Z$.}
\end{align}
Finally, for any $\Lambda$ satisfying \eqref{thm:dioccapf1},
we have in fact
\begin{align*}
\#(A_T\cap\Lambda)=\alpha_F^d(A_T)+O\bigl(\alpha_F^d(A_T)^{1/2+\ve}\bigr)\qquad
\text{as }\: T\to\infty \text{ through $\R$,}
\end{align*}
since for any $n\in\Z^+$ and any $n\leq T\leq n+1$,
we have
$\#(A_n\cap\Lambda)\leq \#(A_T\cap\Lambda)\leq\#(A_{n+1}\cap\Lambda)$ 
and $n=\alpha_F^d(A_n)\leq \alpha_F^d(A_T)\leq\alpha_F^d(A_{n+1})=n+1$.
\end{proof}

Our next purpose is to reduce the above theorem to $\GL_d^1(\mathbb{C})$, where
\[
    \GL_{d}^{1}(\mathbb{C}) = \{g \in \GL_{d}(\mathbb{C}) : \norm{\det(g)}_{\infty} = 1\}.
\]
Note that $\GL_d^1(\CC)$ is unimodular, i.e.\ any left Haar measure of $\GL_d^1(\CC)$ is also a right Haar measure.
We realize $\GL_{d}^{1}(\mathbb{C})$ as a subgroup of $G_d$ via the usual embeddings 
$\GL_{d}^{1}(\mathbb{C}) \hookrightarrow \GL_{d}(\mathbb{C}) \hookrightarrow G_d$.
\begin{lemma}[Reduction to $\GL_{d}^{1}(\mathbb{C})$]\label{lem:reduction}
    With notation as in \autoref{thm:diocca}, 
the relation \eqref{eq:schasymp} holds for (Haar-)almost every $g \in \GL_{d}^{1}(\mathbb{C})$, i.e., for all $\varepsilon > 0$ and for 
almost every $g \in \GL_{d}^{1}(\mathbb{C})$,
    \begin{equation}\label{eq:schasympgl}
        \#\left(A_T \cap gF^{d}\right) = \alpha_F^d(A_{T}) + O{\left(\alpha_F^d(A_{T})^{1/2 + \varepsilon}\right)}
\qquad\text{as }\: T\to\infty.
    \end{equation}
\end{lemma}
\begin{proof}
    For each prime $\mathfrak{p}$ and for our fixed $\bm{v}$, set
\begin{align*}
K'_{\frakp}=\{k\in\GL_d(\scrO_{\frakp})\col k\vecv\in \vecv + \mathfrak{p}^{\alpha_\mathfrak{p}}\mathcal{O}^{d}_\mathfrak{p}\}.
\end{align*}
This is an open subgroup of $\GL_{d}(\mathcal{O}_\mathfrak{p})$.
    Note that for $\mathfrak{p}$ such that $\ord_\mathfrak{p}(\mathfrak{I}) = 0$, i.e., for almost every $\mathfrak{p}$, $K'_\mathfrak{p} = \GL_{d}(\mathcal{O}_\mathfrak{p})$.
    Define $K' = \prod_{\mathfrak{p} < \infty} K'_\mathfrak{p}$.
Then $\GL_{d}(\mathbb{C})K'$ is an open subgroup of $\GL_{d}(\mathbb{A}_F)$,
and $\GL_{d}(\mathbb{C})K' \cap G_d = \GL_{d}^{1}(\mathbb{C}) K'$;
hence $\GL_{d}^{1}(\mathbb{C}) K'$ is an open subgroup of $G_d$.
Therefore, by \autoref{thm:diocca}, 
for $\mu_d^1$-almost every $g \in \GL_{d}^{1}(\mathbb{C}) K'$ we have:
    \[
        \#\left(A_T \cap gF^{d}\right) = \alpha_F^d(A_{T}) + O{\left(\alpha_F^d(A_{T})^{1/2 + \varepsilon}\right)}
\qquad\text{as }\: T\to\infty.
    \]
    Furthermore, since $\GL_{d}^{1}(\mathbb{C}) K' \cong \GL_{d}^{1}(\mathbb{C}) \times K'$, by Fubini's theorem there exists a full measure subset $G' \subseteq \GL_{d}^{1}(\mathbb{C})$ such that for all $g \in G'$ there exists a full measure subset $K'_g \subseteq K'$ so that for all $k \in K'_g$
    \[
        \#\left(A_T \cap gkF^{d}\right) = \alpha_F^d(A_{T}) + O{\left(\alpha_F^d(A_{T})^{1/2 + \varepsilon}\right)}
\qquad\text{as }\: T\to\infty.
    \]

But for any $g\in\GL_d^1(\CC)$, $k\in K'$ and $\vecw\in F^d$,
$gk\vecw\in A_T$ holds if and only if
$g\vecw\in A_{T,\infty}$ and
$k_{\frakp}\vecw\in \vecv+\frakp^{\alpha_{\frakp}}\scrO_{\frakp}^d$
for all $\frakp<\infty$;
and here, since $k_{\frakp}\in K_{\frakp}'$,
the condition
$k_{\frakp}\vecw\in \vecv+\frakp^{\alpha_{\frakp}}\scrO_{\frakp}^d$
is equivalent with
$\vecw\in \vecv+\frakp^{\alpha_{\frakp}}\scrO_{\frakp}^d$.
This shows that 
for any $g\in\GL_d^1(\CC)$, $k\in K'$ and $\vecw\in F^d$,
the two conditions 
$gk\vecw\in A_T$
and $g\vecw\in A_T$ are equivalent,
and therefore $\#(A_T\cap gkF^d)=\#(A_T\cap gF^d)$.
Hence from the above statement, we conclude that for every $g\in G'$,
we have
    \[
        \#\left(A_T \cap gF^{d}\right) = \alpha_F^d(A_{T}) + O{\left(\alpha_F^d(A_{T})^{1/2 + \varepsilon}\right)}
\qquad\text{as }\: T\to\infty.
    \]
\end{proof}

\subsection{Decomposition of the Haar measure of $\GL_{d}^{1}(\mathbb{C})$}\label{ssec:decomphu}
The proof of \autoref{thm:cdiocnew} will follow if, in the case $A_T=E_T$, we can 
reduce \autoref{lem:reduction} further from $\GL_d^1(\C)$ to $U$, where
\[
    U := {\left\{u(\vartheta) : \vartheta \in \M_{m \times n}(\mathbb{C})\right\}}.
\]
Let $H \subseteq \GL_{d}^{1}(\mathbb{C})$ be the following subgroup
\[
    H := {\left\{h = \begin{pmatrix}
        \alpha & 0 \\
        \beta & \gamma
\end{pmatrix} : \alpha \in \GL_{m}(\mathbb{C}), \beta \in \M_{n \times m}(\mathbb{C}), \gamma \in \GL_{n}(\mathbb{C}),
\: \|\det(\alpha)\det(\gamma)\|_\infty=1\right\}}.
\]
We have $H\cap U=\{1_d\}$, and
the multiplication map $\langle h,u\rangle\mapsto hu$ from $H\times U$ to
$\GL_d^1(\CC)$ is a diffeomorphism onto the open subset $HU$ of
$\GL_d^1(\CC)$, the complement of which has Haar measure zero.
Hence by \cite[Theorem 8.32]{aK2002},
parametrizing $\GL_d^1(\CC)$ by $H\times\M_{m\times n}(\CC)$ 
through
$\langle h,\vartheta\rangle\mapsto h\, u(\vartheta)$,
a Haar measure on $\GL_{d}^{1}(\mathbb{C})$ can be expressed 
as the product measure $\d{h}\d{\vartheta}$
where $\d{h}$ is a left Haar measure on $H$ and $\d{\vartheta}$ is the standard Lebesgue measure on
$\M_{m\times n}(\CC)$.
Hence, by Fubini's theorem, \autoref{lem:reduction} implies the following:
\begin{lemma}\label{lem:schasymphue}
    With $A_T$ as in \autoref{thm:diocca}, for all $\varepsilon > 0$ and for almost every $h \in H$,
there exists a full measure subset 
$M_h\subseteq\M_{m\times n}(\CC)$
such that for all $\vartheta\in M_h$,
    \begin{equation}\label{eq:schasymphue}
        \#\left(A_T \cap hu(\vartheta)F^{d}\right) 
= \alpha_F^d(A_{T}) + O{\left(\alpha_F^d(A_{T})^{1/2 + \varepsilon}\right)}
\qquad\text{as }\: T\to\infty.
    \end{equation}
\end{lemma}

\subsection{Proof of \autoref{thm:cdiocnew}}\label{ssec:proofofthmcdiocnew}
In the reduction of \autoref{lem:reduction} to $U$ the following lemma is a key step.
The proof of this lemma is a routine generalization of \cite[Lemma 3.1]{alamghoshyu21}.
\begin{lemma}\label{lem:squeeze}
    There exists $0 < c_0 < \frac{1}{2}$ such that for all $\varepsilon \in (0, c_0)$ there exists an open neighborhood $H_\varepsilon \subseteq H$ of the identity element such that for all $h \in H_\varepsilon$
    \[
        E_{T, \varepsilon}^- \subseteq hE_T \subseteq E_{T, \varepsilon}^+, \qquad \forall T > 10,
    \]
    where
    \[
        E_{T, \varepsilon}^- := E_{T, \varepsilon, \infty}^- \times \prod_{\mathfrak{p} < \infty}(\bm{v} + \mathfrak{p}^{\alpha_\mathfrak{p}}\mathcal{O}_\mathfrak{p}^d)
    \]
    with
    \[
        E_{T, \varepsilon, \infty}^- := {\left\{(\bm{x}, \bm{y}) \in \mathbb{C}^{m} \times \mathbb{C}^{n} : \norm{\bm{x}}_\infty^m \leq (1 + \varepsilon)^{-1} \psi((1 + \varepsilon)\norm{\bm{y}}_\infty^n), \frac{3}{2} \leq \norm{\bm{y}}_\infty^n < (1 + \varepsilon)^{-1}T\right\}};
    \]
    and
    \[
        E_{T, \varepsilon}^+ := E_{T, \varepsilon, \infty}^+ \times \prod_{\mathfrak{p} < \infty}(\bm{v} + \mathfrak{p}^{\alpha_\mathfrak{p}}\mathcal{O}_\mathfrak{p}^d)
    \]
    with
    \begin{gather*}
        E_{T, \varepsilon, \infty}^+ := E_{T, \varepsilon, \infty}' \cup C_0, \\
        E_{T, \varepsilon, \infty}' := {\left\{(\bm{x}, \bm{y}) \in \mathbb{C}^{m} \times \mathbb{C}^{n} : \norm{\bm{x}}_\infty^m \leq (1 + \varepsilon) \psi((1 + \varepsilon)^{-1}\norm{\bm{y}}_\infty^n), \frac{3}{2} \leq \norm{\bm{y}}_\infty^n \leq (1 + \varepsilon)T\right\}}
    \end{gather*}
    and
    \[
        C_0 := {\left\{(\bm{x}, \bm{y}) \in \mathbb{C}^{m} \times \mathbb{C}^{n} : \norm{\bm{x}}_\infty^m \leq 2\psi(1), \frac{1}{2} \leq \norm{\bm{y}}_\infty^n \leq \frac{3}{2}\right\}}.
    \]
\end{lemma}
\begin{remark}
For each $\varepsilon \in (0, c_0)$, both $\{E_{T, \varepsilon, \infty}^{+}\}_{T > 1}$
and $\{E_{T, \varepsilon, \infty}^{-}\}_{T > 1}$ are increasing families.
Furthermore, using similar ideas as in the volume computation of $E_T$ in \S\ref{ssec:dcx},
we see that for all $\varepsilon \in (0, c_0)$ and $T > 10$,
            \begin{equation}\label{eq:voletvep}
                \begin{split}
                    \alpha_F^d(E_{T, \varepsilon}^\pm) &= 
2^d |\Delta_F|^{-d/2} N(\mathfrak{I})^{-d} \pi^d(1 + \varepsilon)^{\pm 1} \int_{\frac{3}{2}}^{(1 + \varepsilon)^{\pm 1} T} \psi((1 + \varepsilon)^{\mp 1}t) \d{t} + O(1) \\
&=2^d |\Delta_F|^{-d/2} N(\mathfrak{I})^{-d} \pi^d(1 + \varepsilon)^{\pm 2} \int_{1}^{T} \psi(t) \d{t} + O(1) \\
                                                 &= (1 + \varepsilon)^{\pm 2} \alpha_F^d(E_T) + O(1).
                \end{split}
            \end{equation}
\end{remark}

\begin{proof}[Proof of \autoref{thm:cdiocnew}]
    Under our assumption in \autoref{thm:cdiocnew} that $\int_{1}^{\infty} \psi(t) \d{t} = \infty$,
we have for all $\varepsilon \in (0, c_0)$:
    \begin{equation}\label{eq:voletvepinfty}
        \alpha_F^d(E_T)\:\text{ and }\: \alpha_F^d(E_{T, \varepsilon}^{\pm}) \to \infty \qquad \text{as} \ T \to \infty.
    \end{equation}
    Let us pick a sequence $\{\varepsilon_\ell\}_{\ell \in \mathbb{N}}$ in $(0,c_0)$ with $\varepsilon_\ell \to0$ as $\ell\to\infty$.
    For each $\ell \in \mathbb{N}$, using
\eqref{eq:voletvepinfty}, and applying \autoref{lem:schasymphue} to $E_{T, \varepsilon_\ell}^\pm$,
we find $h_\ell \in H_{\ve_\ell}$
(where $H_{\ve_\ell}$ is the open set from Lemma~\ref{lem:squeeze}),
and a full measure subset $M_\ell \subseteq \M_{m\times n}(\CC)$ such that for all $\vartheta \in M_\ell$,
    \begin{equation}\label{eq:schasymphuepm}
        \#\left(E_{T, \varepsilon_\ell}^\pm \cap h_\ell u(\vartheta)F^{d}\right) \thicksim \alpha_F^d(E_{T, \varepsilon_\ell}^{\pm}), \qquad \text{as} \ T \to \infty.
    \end{equation}
Furthermore, since $h_\ell \in H_{\ve_\ell}$, we have
    \begin{equation}\label{eq:conleml}
        E_{T, \varepsilon_\ell}^- \subseteq h_\ell E_T \subseteq E_{T, \varepsilon_\ell}^+, \qquad \forall T > 10.
    \end{equation}
    Let $M_\infty := \cap_{\ell \in \mathbb{N}} M_\ell$;
this is still a full measure subset of $\M_{m\times n}(\CC)$.
    For each $\vartheta \in M_\infty$, using \eqref{eq:schasymphuepm}, \eqref{eq:conleml}
and \eqref{eq:voletvep}, 
we see that $\#(E_T \cap u(\vartheta)F^d)$ satisfies the following bounds:
    \[
        \liminf_{T \to \infty} \frac{\#(E_T \cap u(\vartheta)F^d)}{\alpha_F^d(E_T)} \geq \liminf_{T \to \infty} \frac{\#\left(E_{T, \varepsilon_\ell}^- \cap h_\ell u(\vartheta)F^{d}\right)}{\alpha_F^d(E_T)} = (1 + \varepsilon_\ell)^{-2},
    \]
    and
    \[
        \limsup_{T \to \infty} \frac{\#(E_T \cap u(\vartheta)F^d)}{\alpha_F^d(E_T)} \leq \limsup_{T \to \infty} \frac{\#\left(E_{T, \varepsilon_\ell}^+ \cap h_\ell u(\vartheta) F^{d}\right)}{\alpha_F^d(E_T)} = (1 + \varepsilon_\ell)^{2}.
    \]
    Letting $\ell \to \infty$ we find our desired result.
\end{proof}




\bibliographystyle{alpha}
\bibliography{bib}

\begin{thebibliography}{GKY20}

\bibitem[AG20]{alamghosh20}
Mahbub Alam and Anish Ghosh.
\newblock {Equidistribution on homogeneous spaces and the distribution of
  approximates in {D}iophantine approximation}.
\newblock {\em Trans. Amer. Math. Soc.}, 373(5):3357--3374, 2020.

\bibitem[AGH23]{alamghoshhan21}
M.~Alam, A.~Ghosh, and J.~Han.
\newblock {Higher moment formulae and limiting distributions of lattice
  points}.
\newblock {\em To appear in the J. of the Inst. of Math. of Jussieu}, 2023.
\newblock preprint, \url{https://arxiv.org/abs/2111.00848}.

\bibitem[AGY21]{alamghoshyu21}
Mahbub Alam, Anish Ghosh, and Shucheng Yu.
\newblock Quantitative {D}iophantine approximation with congruence conditions.
\newblock {\em J. Th\'{e}or. Nombres Bordeaux}, 33(1):261--271, 2021.

\bibitem[AM09]{jathreyagmargulis09}
Jayadev~S. Athreya and Gregory~A. Margulis.
\newblock Logarithm laws for unipotent flows. {I}.
\newblock {\em J. Mod. Dyn.}, 3(3):359--378, 2009.

\bibitem[BG19]{bjorklundgorodnik19}
Michael Björklund and Alexander Gorodnik.
\newblock {Central limit theorems for Diophantine approximants}.
\newblock {\em Math. Annalen}, 374:1371--1437, 2019.

\bibitem[Bor91]{aB91}
Armand Borel.
\newblock {\em Linear algebraic groups}.
\newblock Springer-Verlag, New York, second edition, 1991.

\bibitem[For18]{ford18}
Lester~R. Ford.
\newblock {Rational Approximations to Irrational Complex Number}.
\newblock {\em Transactions of the American Mathematical Society}, 19(1):1--42,
  1918.

\bibitem[For25]{ford25}
Lester~R. Ford.
\newblock {On the closeness of approach of complex rational fractions to a
  complex irrational number}.
\newblock {\em Trans. Amer. Math. Soc.}, 27(2):146--154, 1925.

\bibitem[Gar23]{nG2023c}
N.~Gargava.
\newblock {\em Mean Value Theorems for Collections of Lattices with a
  Prescribed Group of Symmetries}.
\newblock PhD thesis, EPFL, 2023.

\bibitem[GKY20]{ghoshkelmeryu2020}
A.~Ghosh, D.~Kelmer, and S.~Yu.
\newblock {Effective Density for Inhomogeneous Quadratic Forms I: Generic Forms
  and Fixed Shifts}.
\newblock {\em Int. Math. Res. Not. IMRN}, 08 2020.
\newblock rnaa206.

\bibitem[Gro38]{groshev38}
A.~V. Groshev.
\newblock {A theorem on system of linear forms}.
\newblock {\em Doklady Akad. Nauk SSSR}, 19:151--152, 1938.
\newblock in Russian.

\bibitem[GSV23]{gargavasv23}
N.~Gargava, V.~Serban, and M.~Viazovska.
\newblock {Moments of the number of points in a bounded set for number field
  lattices}.
\newblock arXiv:2308.15275, 2023.

\bibitem[Han22]{han22}
Jiyoung Han.
\newblock {Rogers' mean value theorem for {$S$}-arithmetic {S}iegel transforms
  and applications to the geometry of numbers}.
\newblock {\em J. Number Theory}, 240:74--106, 2022.

\bibitem[Hug23]{nH2023}
N.~Hughes.
\newblock Mean values over lattices in number fields and effective diophantine
  approximation.
\newblock arXiv:2306.02499, 2023.

\bibitem[Hur87]{hurwitz1887}
A.~Hurwitz.
\newblock {\"{U}ber die {E}ntwicklung complexer {G}r\"{o}ssen in
  {K}ettenbr\"{u}che}.
\newblock {\em Acta Math.}, 11(1-4):187--200, 1887.

\bibitem[Khi26]{khinchin26}
A.~Khinchin.
\newblock {Zur metrischen theorie der diophantischen approximationen}.
\newblock {\em Math. Z.}, 24:706--14, 1926.

\bibitem[Kim24]{kim24}
Seungki Kim.
\newblock {Adelic Rogers' integral formula}.
\newblock {\em Journal of the London Mathematical Society}, 109(1):e12830,
  2024.

\bibitem[KL16]{kleinbockly16}
Dmitry Kleinbock and Tue Ly.
\newblock {Badly approximable {$S$}-numbers and absolute {S}chmidt games}.
\newblock {\em J. Number Theory}, 164:13--42, 2016.

\bibitem[Kna02]{aK2002}
A.~W. Knapp.
\newblock {\em Lie groups beyond an introduction}, volume 140 of {\em Progress
  in Mathematics}.
\newblock Birkh\"auser Boston Inc., Boston, MA, second edition, 2002.

\bibitem[KST17]{kleinbockshitomanov17}
D.~Kleinbock, R.~Shi, and G.~Tomanov.
\newblock {$S$-adic version of Minkowski's geometry of numbers and Mahler's
  compactness criterion}.
\newblock {\em J. Number Theory}, 174:150--163, 2017.

\bibitem[KY19]{dKsY2019}
D.~Kelmer and S.~Yu.
\newblock The second moment of the siegel transform in the space of symplectic
  lattices.
\newblock {\em International Mathematics Research Notices}, 02 2019.

\bibitem[LeV52]{leveque52}
W.~J. LeVeque.
\newblock {Continued fractions and approximations in {$k(i)$}. {I}, {II}}.
\newblock {\em Indag. Math.}, 14:526--535, 536--545, 1952.
\newblock Nederl. Akad. Wetensch. Proc. Ser. A {\bf 55}.

\bibitem[Ly16]{ly16}
Tue~Ngoc Ly.
\newblock {\em Diophantine {A}pproximation in {A}lgebraic {N}umber {F}ields and
  {F}lows on {H}omogeneous {S}paces}.
\newblock PhD thesis, Brandeis University, 2016.
\newblock Thesis (Ph.D.).

\bibitem[MR58]{macbeathrogers58}
A.~M. Macbeath and C.~A. Rogers.
\newblock {Siegel's Mean Value Theorem in the Geometry of Numbers}.
\newblock {\em Math. Proc. Cam. Phil. Soc.}, 54(2):139--151, 1958.

\bibitem[MS10]{MarklofStrombergsson2010}
Jens Marklof and Andreas Str{\"o}mbergsson.
\newblock {The distribution of free path lengths in the periodic {L}orentz gas
  and related lattice point problems}.
\newblock {\em Ann. of Math. (2)}, 172(3):1949--2033, 2010.

\bibitem[Nak88]{nakada88}
Hitoshi Nakada.
\newblock {On metrical theory of Diophantine approximation over imaginary
  quadratic field}.
\newblock {\em Acta Arith.}, 51:399--403, 1988.

\bibitem[NRS20]{nesharimruhrshi2020}
E.~Nesharim, R.~R\"{u}hr, and R.~Shi.
\newblock {Metric {D}iophantine approximation with congruence conditions}.
\newblock {\em Int. J. Number Theory}, 16(9):1923--1933, 2020.

\bibitem[Oes84]{jO84}
Joseph Oesterl\'{e}.
\newblock Nombres de tamagawa et groupes unipotents en caract\'{e}ristique
  {$p$}.
\newblock {\em Invent. Math.}, 78(1):13--88, 1984.

\bibitem[Qu{\^e}91]{rqueme89}
Roland Qu{\^e}me.
\newblock {On {D}iophantine approximation by algebraic numbers of a given
  number field: a new generalization of {D}irichlet approximation theorem}.
\newblock Number 198-200, pages 273--283. 1991.
\newblock Journ\'{e}es Arithm\'{e}tiques, 1989 (Luminy, 1989).

\bibitem[Rog55]{rogers55a}
C.~A. Rogers.
\newblock {Mean values over the space of lattices}.
\newblock {\em Acta Mathematica}, 94:249--287, 1955.

\bibitem[Rog56]{rogers56}
C.~A. Rogers.
\newblock {The Number of Lattice Points in a Set}.
\newblock {\em Proc. of the Lon. Math. Soc.}, s3-6:305--320, 1956.

\bibitem[Rud87]{rudin87}
W.~Rudin.
\newblock {\em {Real and Complex Analysis}}.
\newblock McGraw-Hill, 1987.

\bibitem[Sch57]{wS57}
W.~M. Schmidt.
\newblock Mittelwerte \"uber {G}itter.
\newblock {\em Monatsh. Math.}, 61:269--276, 1957.

\bibitem[Sch58a]{schmidt58b}
W.~Schmidt.
\newblock {The Measure of the Set of Admissible Lattices}.
\newblock {\em Proceedings of the American Mathematical Society},
  9(3):390--403, 1958.

\bibitem[Sch58b]{wS58a}
W.~M. Schmidt.
\newblock Mittelwerte \"uber {G}itter. {II}.
\newblock {\em Monatsh. Math.}, 62:250--258, 1958.

\bibitem[Sch58c]{schmidt58a}
Wolfgang Schmidt.
\newblock {On the convergence of mean values over lattices}.
\newblock {\em Canadian J. Math.}, 10:103--110, 1958.

\bibitem[Sch60a]{schmidt60b}
W.~Schmidt.
\newblock {A metrical theorem in Diophantine approximation}.
\newblock {\em Can. J. Math.}, 12:619--31, 1960.

\bibitem[Sch60b]{schmidt60a}
W.~M. Schmidt.
\newblock {A metrical theorem in geometry of numbers}.
\newblock {\em Trans. Amer. Math. Soc.}, 95:516--529, 1960.

\bibitem[Sch67]{schmidt67}
Wolfgang~M. Schmidt.
\newblock {On heights of algebraic subspaces and diophantine approximations}.
\newblock {\em Ann. of Math. (2)}, 85:430--472, 1967.

\bibitem[Sch75]{aschmidt75}
Asmus~L. Schmidt.
\newblock {Diophantine approximation of complex numbers}.
\newblock {\em Acta Math.}, 134:1--85, 1975.

\bibitem[Sie45]{siegel45}
C.~L. Siegel.
\newblock {A Mean Value Theorem in Geometry of Numbers}.
\newblock {\em Ann. Math. Princeton}, 46:340--347, 1945.

\bibitem[S{\"o}d11]{aS2011o}
A.~S{\"o}dergren.
\newblock On the {P}oisson distribution of lengths of lattice vectors in a
  random lattice.
\newblock {\em Math. Z.}, 269(3-4):945--954, 2011.

\bibitem[Spr79]{vS79}
V.~G. Sprindzuk.
\newblock {\em Metric theory of Diophantine approximations}.
\newblock V. H. Winston \& Sons, Washington, D.C., 1979.
\newblock Translated from the Russian and edited by Richard A. Silverman, With
  a foreword by Donald J. Newman, Scripta Series in Mathematics.

\bibitem[SS06]{pSaS2006}
P.~Sarnak and A.~Str\"ombergsson.
\newblock Minima of epstein's zeta function and heights of flat tori.
\newblock {\em Invent. Math.}, 165:115--151, 2006.

\bibitem[SS19]{strombergssonsodergren19}
Andreas Strömbergsson and Anders Södergren.
\newblock {On the generalized circle problem for a random lattice in large
  dimension}.
\newblock {\em Adv. in Math.}, 345:1042--1074, 2019.

\bibitem[SS22]{aSaS2022}
A.~Strömbergsson and A.~Södergren.
\newblock On a mean value formula for multiple sums over a lattice and its
  dual.
\newblock arXiv:2211.05454, 2022.

\bibitem[Sul82]{sullivan82}
Dennis Sullivan.
\newblock {Disjoint spheres, approximation by imaginary quadratic numbers, and
  the logarithm law for geodesics}.
\newblock {\em Acta Math.}, 149:215--237, 1982.

\bibitem[Sz{\"u}62]{szusz62}
P.~Sz{\"u}sz.
\newblock {\"{U}ber die metrische {T}heorie der diophantischen {A}pproximation.
  {II}}.
\newblock {\em Acta Arith.}, 8:225--241, 1962.

\bibitem[Wei82]{weil82}
André Weil.
\newblock {\em {Adeles and algebraic groups}}.
\newblock Birkhäuser Boston, Mass., 1982.
\newblock With appendices by M.\ Demazure and Takashi Ono.

\bibitem[Wei95]{weil95}
Andre Weil.
\newblock {\em {Basic Number Theory}}.
\newblock Classics in Mathematics. Springer, 1995.

\end{thebibliography}

\end{document}